\documentclass[leqno,11pt]{article}
\usepackage{latexsym}
\usepackage{amssymb}
\usepackage{amsfonts}
\usepackage{amsmath}
\usepackage{epsfig}
\usepackage{color}

\allowdisplaybreaks
\makeatletter
\long\def\unmarkedfootnote#1{{\long\def\@makefntext##1{##1}\footnotetext{#1}}}
\makeatother

\newtheorem{definition}{Definition}[section]

\newtheorem{lemma}[definition]{Lemma}

\newtheorem{theorem}[definition]{Theorem}

\newtheorem{proposition}[definition]{Proposition}

\newtheorem{remark}[definition]{Remark}

\def\o{\Omega}

\def\m2{|\Omega | /2}
\def\M2{\frac{|\Omega |}{2}}
\def\u+{u_+^*}

\def\-p{\overline{p}}

\def\w0{{W_0^{1,p}(\Omega)}}

\def\MR{\mathcal R}
\def\R{\mathbb R}
\def\N{\mathbb N}

\def\ep{\varepsilon}

\def\rn{{{\R}^n}}
\def\T{T_{t}}

\newcommand{\hh}{{\cal H}^{n-1}}
\newcommand{\medint}{-\kern  -,395cm\int}
\newcommand{\medintinrigo}{-\kern  -,315cm\int}
\newcommand{\medelle}{-\kern  -,235cm L}
\newcommand{\medellenrigo}{-\kern  -,180cm L}
\newcommand{\qed}{\thinspace\null\nobreak\hfill
\hbox{\vbox{\kern-.2pt\hrule height.2pt
depth.2pt\kern-.2pt\kern-.2pt \hbox to1.8mm {\kern-.2pt\vrule
width.4pt \kern-.2pt\raise1.8mm\vbox to.2pt{} \lower0pt\vtop
to.2pt{}\hfil\kern-.2pt \vrule
width.4pt\kern-.2pt}\kern-.2pt\kern-.2pt \hrule height.2pt
depth.2pt \kern-.2pt}}\par\medbreak}

\title{Second-order $L^2$-regularity in nonlinear elliptic problems
} 
\numberwithin{equation}{section}

\author{
  Andrea Cianchi\\
 {\it Dipartimento di Matematica e Informatica \lq\lq U. Dini", Universit\`a di Firenze}\\ {\it Viale Morgagni 67/A, 50134 Firenze, Italy} \\{\it  e-mail: cianchi@unifi.it}
\bigskip
\\
  Vladimir G. Maz'ya \\
  {\it   Department of Mathematics, Link\"oping University, SE-581
83 Link\"oping, Sweden}
  \\ and \\
 {\it  Department of Mathematical Sciences, M\&O Building}\\ {\it University of Liverpool, Liverpool L69 3BX,
 UK;}
\\ {\it e-mail: vlmaz@mai.liu.se}
}

\date{}

\begin{document}
\maketitle
\begin{abstract}
A second-order regularity theory is developed for solutions to a class of quasilinear elliptic equations in divergence form, including the $p$-Laplace equation, with merely square-integrable right-hand side. Our results amount to the existence and square integrability of the weak derivatives of the nonlinear expression of the gradient  under the divergence operator. This provides  a  nonlinear counterpart of  the classical $L^2$-coercivity theory  for linear problems,
which is  missing in the existing literature. 
Both local  and global estimates are established. The latter apply to solutions to either Dirichlet or Neumann boundary value problems. Minimal regularity  on the boundary of the domain is required. If the domain is convex, no regularity of its boundary is needed  at all.
\end{abstract}

\unmarkedfootnote {

\smallskip
\par\noindent 
This research was partly supported by the Research Project of the
Italian Ministry of University and Research (MIUR) Prin 2012 n.2012TC7588 "Elliptic and
parabolic partial differential equations: geometric aspects, related
inequalities, and applications", and by  GNAMPA   of INdAM (National Institute of High Mathematics).}

\section{Introduction}\label{sec1}

A  prototypal
result in the theory of elliptic equations asserts that, if $\Omega$ is a bounded open set in $\rn $, $n \geq 2$,  with $\partial \Omega \in C^2$, and $u$ is the weak solution to the Dirichlet problem for the inhomegenous Laplace equation whose right-hand side  $f \in L^2(\Omega)$, then  $u\in W^{2,2}(\Omega)$. Moreover, a two-sided  coercivity estimate for $\|\nabla ^2u\|_{L^2(\Omega)}$ holds in terms of  $\|f\|_{L^2(\Omega)}$, up to multiplicative constants.
%
%
This can be traced back to \cite{Bernstein} for $n=2$, and to \cite{Schauder} for $n \geq 3$. A   comprehensive   analysis   of this topic can be found  in  \cite{ADN}, \cite[Chapter 10]{Hor}, \cite[Chapter 3]{LU}, \cite[Chapter 14]{MaSh}.
\par The regularity theory  for (possibly degenerate or singular) nonlinear equations in divergence form, extending the  Laplace equation, whose model is the  $p$-Laplace equation, has thoroughly been developed in the last fifty years. 
Regularity properties of solutions and of their first-order derivatives have been investigated in a number of contributions, including the classics \cite{CDiB, Di,  DiMa, Ev, Iw, KM,  Le, Li, Simon, To, Uhl, Ur} and the more recents advances \cite{BCDKS, BDS, CKP, cmARMA, BDS, DM,   KuM}.
\par 
Despite the huge amount of work devoted to this kind of equations,
the picture of second-order regularity for their solutions 
is apparently still  quite incomplete. A result is available for $p$-harmonic functions, namely local solutions $u$ to the homogenous equation
$$- {\rm div} (|\nabla u|^{p-2}\nabla u ) = 0 \quad {\rm in}\,\,\, \o\,,$$
and asserts that the nonlinear expression of the gradient $|\nabla u|^{\frac{p-2}2}\nabla u \in W^{1,2}_{\rm loc}(\Omega)$ -- see \cite{Uhl} for $p \in (2, \infty)$, and \cite{CDiB} for every $p \in (1, \infty)$. If $p\in (1,2)$, coupling this property with the local boundedness  of $\nabla u$   in $\o$ ensures that $u \in W^{2,2}_{\rm loc}(\Omega)$. On the other hand, the existence of second-order weak  derivatives of $p$-harmonic functions is an open problem for $p \in (2,\infty)$. 
\par Information on this issue concerning inhomogeneous equations is  even more limited. In fact, this case   seems to be almost unexplored.
With this regard, let us mention that  (global) twice weak differentiability of solutions to Dirichlet problems for the inhomogeneous  $p$-Laplace equation is proved in  \cite{BC}  under the assumption that  $p$ is smaller than, and sufficiently close to $2$,
%
and relies upon the linear theory, via a perturbation argument. Fractional-order regularity  for the gradient of solutions to a class of nonlinear inhomogeneous equations, modelled upon the  $p$-Laplacian, is established in \cite{Mi1}. An earlier contribution in this direction is \cite{SimonJ}
\par
The present paper offers
a second-order  regularity principle
for a class of quasilinear elliptic problems in divergence form,   that encompasses  the inhomegenous $p$-Laplace equation 
$$- {\rm div} (|\nabla u|^{p-2}\nabla u ) = f(x) \quad {\rm in}\,\,\, \o\,,$$
 for any $p \in (1, \infty)$ and any right-hand side   $f\in L^2(\o)$.  
In contrast with   the customary results recalled above, our statements   involve exactly the    nonlinear function of $\nabla u$ appearing under divergence in the relevant elliptic operators. In the light of our conclusions, this turns out to be   the correct expression to call into play, inasmuch as it admits a two-sided $L^2$-estimate in terms of the datum on the right-hand side of the equation, and hence 
%
%
exhibits a  regularity-preserving property. 
%
%
%
  \par Both local solutions, and solutions to Dirichlet and  Neumann boundary value problems are addressed.
A distinctive trait of our results is the minimal regularity imposed on $\partial \Omega$ when dealing with global bounds. In particular, if $\Omega$ is convex, no additional regularity has to be required on $\partial \Omega$.  However, we stress that  the results to be proved are new even for smooth domains. 
\par An additional striking feature is that they apply to a very weak notion of solutions, which has  to be adopted since the right-hand side of the equations is allowed to enjoy a low degree of integrability.
%
\par To conclude this preliminary overview, let us point out that the validity of second-order $L^2$-estimates raises the natural question of a more general second-order theory in $L^q$  for $q \neq 2$, or in other function spaces. 

\section{Main results}\label{results}

Although our main focus is on global estimates for solutions to boundary value problems,  
we begin our discussion with a local bound  for local solutions, of independent interest.  The equations under consideration   have the form
%
\begin{equation}\label{localeq}
- {\rm div} (a(|\nabla u|)\nabla u ) = f(x)  \quad {\rm in}\,\,\, \o \\
\end{equation}
where $\Omega$ is any open set in $\mathbb R^n$, and $f \in L^2_{\rm loc}(\Omega)$. The function $a : (0, \infty ) \to (0, \infty )$ is of class  $C^1(0, \infty )$, and such that
\begin{equation}\label{infsup}
-1 < i_a \leq s_a < \infty,
\end{equation}
where
 \begin{equation}\label{ia}
i_a= \inf _{t >0} \frac{t a'(t)}{a(t)} \qquad \hbox{and} \qquad
s_a= \sup _{t >0} \frac{t a'(t)}{a(t)},
\end{equation}
and $a'$ stands for the derivative of $a$.
Assumption \eqref{infsup} ensures that the differential operator in \eqref{localeq} satisfies   ellipticity and monotonicity conditions, not necessarily of power type \cite{cmCPDE, cmARMA}. Regularity for equations governed by   generalized nonlinearities of this kind has also been extensively studied -- see e.g.  \cite{Baroni, BreitSV, cianchibound, CianchiAIHP, DKS, DieningSV, Kor, Li, Mar, Talenti_orlicz}.
Observe that the standard $p$-Laplace operator corresponds to
the choice $a(t)=t^{p-2}$, with $p >1$. Clearly,   $i_a=s_a=p-2$ in this case.
\par
As already warned in Section \ref{sec1}, due to the mere square summability assumption on the function $f$, solutions to equation \eqref{localeq} may have to be understood in a suitable generalized sense, even in the case of the $p$-Laplacian. We  shall further comment on this at the end of this section. Precise definitions  can be found in Sections \ref{proofs} and \ref{proofsloc}.
\par In what follows, $B_r(x)$ denotes the ball with radius $r>0$, centered at $x \in \rn$. The simplified notation $B_r$ is employed when information on the center is irrelevant. In this case, balls with different radii appearing in the same formula  (or proof) will be tacitly assumed to have the same center.
\begin{theorem}\label{secondloc} {\rm {\bf [Local estimate]}}
Assume that the function $a \in C^1(0, \infty)$, and  satisfies condition \eqref{infsup}. Let $\Omega$ be any open set in $\mathbb R^n$, with $n \geq 2$, and let $f \in L^2_{\rm loc}(\Omega)$.
   Let $u$ be a generalized  local solution   to  equation
\eqref{localeq}. 
Then
\begin{equation}\label{secondloc1} a(|\nabla u|) \nabla u
\in W^{1,2}_{\rm loc}(\Omega ),
\end{equation}
and there exists a constant $C=C(n, i_a, s_a)$ such that
\begin{equation}\label{secondloc2} \|a(|\nabla u|) \nabla
u\|_{W^{1,2}(B_R)} \leq C \big(\|f\|_{L^2(B_{2R})} + R^{-\frac n2}\|a(|\nabla u|) \nabla
u\|_{L^1(B_{2R})} \big)
\end{equation}
for any ball $B_{2R} \subset \subset \Omega$.
\end{theorem}

%

\begin{remark}\label{plaplace}{\rm
Observe that the  expression $a(|\nabla u|)\nabla u$   agrees with $|\nabla u|^{p-2}\nabla u$ when the differential operator in equation \eqref{localeq} is the $p$-Laplacian,  and hence differs in the exponent of $|\nabla u|$ from the results recalled above about $p$-harmonic functions.}
\end{remark}

\smallskip
\par
Our global results concern  Dirichlet or Neumann problems, with homogeneous boundary data, associated with equation \eqref{localeq}. Namely,  Dirichlet problems of the form
\begin{equation}\label{eqdirichlet}
\begin{cases}
- {\rm div} (a(|\nabla u|)\nabla u ) = f(x)  & {\rm in}\,\,\, \o \\
 u =0  &
{\rm on}\,\,\,
\partial \o \,
\end{cases}
\end{equation}
and   Neumann problems of the form
\begin{equation}\label{eqneumann}
\begin{cases}
- {\rm div} (a(|\nabla u|)\nabla u ) = f(x)  & {\rm in}\,\,\, \o \\
 \displaystyle  \frac{\partial u}{\partial \nu} =0  &
{\rm on}\,\,\,
\partial \o \,.
\end{cases}
\end{equation}
 Here, $\Omega$ is a bounded open set in $\rn$,  $\nu$ denotes the outward unit vector on $\partial \Omega$,  $f \in L^2(\Omega)$,
and $a : (0, \infty ) \to (0, \infty )$ is as above. Of course, the compatibility condition 
\begin{equation}\label{0mean}
\int _\o f(x)\, dx =0
\end{equation}
 has to be required when dealing with \eqref{eqneumann}. 
\par A basic version of the  global second-order estimates for the solutions to 
\eqref{eqdirichlet} and \eqref{eqneumann} 
holds in
any bounded convex open set  $\o \subset \rn$. 


\begin{theorem}\label{secondconvex} {\rm {\bf [Global estimate in convex domains]}}
Assume that the function $a \in C^1(0, \infty)$, and satisfies condition \eqref{infsup}. Let $\o$ be any convex bounded open set in $\rn$, with $n \geq 2$, and let    $f \in
L^2(\Omega )$. 
Let
$u$ be the generalized  solution   to  either the
Dirichlet problem \eqref{eqdirichlet}, or   the Neumann problem
\eqref{eqneumann}. 
Then
\begin{equation}\label{secondconvex1} a(|\nabla u|) \nabla u
\in W^{1,2}(\Omega ).
\end{equation}
Moreover,
\begin{equation}\label{secondconvex2} C_1 \|f\|_{L^2(\Omega)} \leq \|a(|\nabla u|) \nabla
u\|_{W^{1,2}(\Omega)} \leq C_2 \|f\|_{L^2(\Omega)}
\end{equation}
for some constants $C_1=C_1(n,  s_a)$ and $C_2=C_2(\o, i_a, s_a)$.
\end{theorem}

\par
Heuristically speaking, the validity of a global estimate in Theorem \ref{secondconvex} is related to the fact that the  second fundamental form on the boundary of a convex set is semidefinite. 
In the main 
 result of this paper, the convexity assumption on $\Omega$ is abandoned. Dropping   signature information  on the  (weak) second fundamental form on $\partial \Omega$ calls for an assumption on its summability.
%
We assume that the domain $\Omega$ is    locally the subgraph of a Lipschitz continuous function of $(n-1)$ variables, which is also twice weakly differentiable. The second-order derivatives of this function are required to belong to either the weak Lebesgue space $L^{n-1}$, called $L^{n-1, \infty}$, or the weak Zygmund space $L\log L$, called  $L^{1, \infty} \log L$, according to whether $n \geq 3$ or $n=2$. This will be denoted by $\partial \Omega \in  L^{n-1, \infty}$, and  $\partial \Omega  \in L^{1, \infty} \log L$, respectively.
 As a consequence,  the weak second fundamental form $\mathcal B$ on $\partial \Omega$
belongs to the same weak type spaces with respect to the $(n-1)$-dimensional Hausdorff measure $\mathcal H^{n-1}$ on $\partial \Omega$ . Our key summability assumption on $\mathcal B$ amounts to:
\begin{equation}\label{smalln}
 \lim _{r\to 0^+} \Big(\sup _{x \in \partial \Omega} \|\mathcal B \|_{L^{n-1, \infty}(\partial \Omega \cap B_r(x))}\Big) < c   \quad \hbox{if $n \geq 3$,}
\end{equation}
or
\begin{equation}\label{small2}
\lim _{r\to 0^+} \Big(\sup _{x \in \partial \Omega} \|\mathcal B\|_{L^{1, \infty} \log L(\partial \Omega \cap B_r(x))}\Big)< c   \quad \hbox{if $n =2$,}
\end{equation}
for a suitable constant $c=c(L_\Omega, d_\Omega, n, i_a, s_a)$. Here,  $L_\o$ denotes the Lipschitz constant  of $\o$, and  $d_\o$ its diameter. Let us emphasize that such an assumption is  essentially sharp -- see Remark \ref{remsharp} below.

\begin{theorem}\label{seconddir} {\rm {\bf [Global estimate in minimally regular domains]}}
Assume that the function $a \in C^1(0, \infty)$, and  satisfies condition \eqref{infsup}. Let $\Omega$ be a Lipschitz bounded domain in $\rn$, $n \geq 2$ such that   $\partial \Omega \in  W^2L^{n-1, \infty}$ if $n \geq 3$, or $\partial \Omega  \in W^2L^{1, \infty} \log L$ if $n =2$ .
  Assume that $f \in
L^2(\Omega )$, and let $u$ be the generalized  solution   to either  the
Dirichlet problem \eqref{eqdirichlet}, or   the Neumann problem
\eqref{eqneumann}. There exists a constant $c=c(L_\Omega, d_\o, n, i_a, s_a)$ such that, if  $\Omega$ fulfils \eqref{smalln} or \eqref{small2} for such a constant $c$,
then
\begin{equation}\label{seconddir1} a(|\nabla u|) \nabla u
\in W^{1,2}(\Omega ).
\end{equation}
Moreover,
\begin{equation}\label{seconddir2} C_1 \|f\|_{L^2(\Omega)} \leq \|a(|\nabla u|) \nabla
u\|_{W^{1,2}(\Omega)} \leq C_2 \|f\|_{L^2(\Omega)}
\end{equation}
for some positive constants $C_1=C_1(n,  s_a)$ and $C_2=C_2(\o, i_a, s_a)$.
\end{theorem}

We conclude this section with some remarks   on Theorems \ref{secondloc}, \ref{secondconvex} and \ref{seconddir}.

\begin{remark}\label{remsharp}
 {\rm
Assumption \eqref{smalln}, or \eqref{small2}, cannot be weakened in Theorem \ref{seconddir} for all equations of the form appearing in \eqref{eqdirichlet} and
\eqref{eqneumann}.  This can be shown by taking into account  the linear problem corresponding to the case when the function $a$ is constant.
 Indeed, domains $\o$  can be exhibited such that $\partial \o \in W^2L^{n-1, \infty}$ if $n \geq 3$  \cite{Ma73}, or $\partial \o \in  W^2L^{1, \infty} \log L$ if $n =2$ \cite{Ma67}, but   the limit in \eqref{smalln} or \eqref{small2}   exceeds some explicit threshold,  and  the corresponding solution $u$ to the Dirichlet problem for the Laplace equation  fails to belong to $W^{2,2}(\o)$ (see also   \cite[Section 14.6.1]{MaSh} in this connection). 
}
\end{remark}

\begin{remark}\label{rembound}
 {\rm  
Condition \eqref{smalln}  is certainly fulfilled if $\partial \o \in W^{2,n-1}$, and \eqref{small2} is fulfilled if $\partial \o \in W^{2}L\log L$, or, a fortiori, if  $\partial \o \in W^{2,q}$ for some $q>1$. This follows from the embedding of $L^{n-1}$ into $L^{n-1, \infty}$ and of $L\log L$ (or $L^q$) into $L^{1, \infty}\log L$ for $q>1$, and from the absolute continuity of the norm in any Lebesgue and Zygmund space.
Notice also that, 
since the Lorentz space $L^{n-1,1} \subsetneqq L^{n-1}$, assumption 
 \eqref{smalln}
is, in particular, weaker than requiring that $\partial \o \in W^{2}L^{n-1,1}$. The latter condition  has been shown to ensure the global boundedness of the gradient of the solutions to problems \eqref{eqdirichlet} or \eqref{eqneumann}, for $n \geq 3$, provided that $f$ belongs to the Lorentz space $L^{n,1}(\o)$ \cite{cmCPDE, cmARMA}. Note that hypothesis \eqref{smalln} does not imply that $\partial \Omega \in C^{1,0}$, a property that is instead certainly fulfilled under the stronger condition that $\partial \o \in W^{2}L^{n-1,1}$.
}
\end{remark}

\begin{remark}\label{remW22}
 {\rm
\par \color{black} The gloal gradient bound mentioned in Remark \ref{rembound} enables one  to show, via a minor variant in the proof of Theorems \ref{secondconvex}--\ref{seconddir}, that the solutions to problems \eqref{eqdirichlet} and \eqref{eqneumann} are actually in $W^{2,2}(\o)$, provided that 
\begin{equation}\label{infa}
\inf _{t\in[0, M]}a(t) > 0
\end{equation}
for every $M>0$, and $f$ and $\o$ have the required regularity for the relevant gradient bound to hold. A parallel result holds for local solutions to the equation \eqref{localeq}, thanks to a local gradient estimate from \cite{Baroni}, extending \cite{DM}. To be more specific, if $f \in L^{n,1}_{\rm loc}(\Omega)$, and
    $u$ is a generalized  local solution   to  equation
\eqref{localeq},
then
\begin{equation}\label{W221}  u
\in W^{2,2}_{\rm loc}(\Omega ).
\end{equation}
Moreover, if $n \geq 3$,  $f \in L^{n,1}(\Omega)$,  $\partial \Omega \in W^2L^{n-1,1}$, and   $u$ is the generalized  solution   to either  the
Dirichlet problem \eqref{eqdirichlet}, or   the Neumann problem
\eqref{eqneumann}, 
then
\begin{equation}\label{W223}  u
\in W^{2,2}(\Omega ).
\end{equation}
Equation \eqref{W223} continues to hold if $\Omega$ is any bounded convex domain in $\rn$, whatever $\partial \Omega$ is.
\\ Let us stress that these conclusion may fail if assumption \eqref{infa} is dropped. This can be verified, for instance, on choosing $a(t)= t^{p-2}$, i.e. the $p$-Laplace operator, and considering functions of the form $u(x)= |x_1|^\beta$, where $x=(x_1, \dots , x_n)$ and $\beta >1$. These functions are local solutions to \eqref{localeq} with $f \in L^{n,1}_{\rm loc}(\rn)$ (and even $f \in L^{\infty}_{\rm loc}(\rn)$)  provided that $p$ is large enough, but $u \notin W^{2,2}_{\rm loc}(\rn)$ if $\beta \leq \tfrac 32$. In fact, $u \notin W^{2,q}_{\rm loc}(\rn)$ for any given $q>1$, if $\beta$ is sufficiently close to $1$.}
\end{remark}

\color{black}

\begin{remark}\label{remsol}
 {\rm 
Weak solutions to problems \eqref{eqdirichlet} or \eqref{eqneumann}, namely  distributional solutions  belonging to the energy space associated with the relevant differential operator, need not exist if $f$ is merely in $L^2(\o)$. This phenomenon is well-known to occur   in the model case of the $p$-Laplace equation, if $p$ is not large enough for $L^2(\Omega)$ to be contained in the dual of $W^{1,p}(\Omega)$. Yet, weaker definitions of solutions to boundary value problems for this equation, ensuring their uniqueness, which  apply to any $p \in (1, \infty)$ and even to right-hand sides $f \in L^1(\o)$,  are available in the literature\cite{ACMM, BBGGPV, BG, DallA, DM, LM, Ma69, Mu}.  Among the diverse, but a posteriori equivalent, definitions, we shall adopt that (adjusted to the  framework under consideration in this paper)  of a solution which is the limit of a sequence of solutions to problems whose right-hand sides are smooth and converge to $f$ \cite{DallA}. This will be called a generalized solution throughout. A parallel notion of  generalized local solution to \eqref{localeq} will be empolyed. 
A generalized solution  need not be  weakly differentiable. However, it is associated with a vector-valued function on $\o$, which plays the role of a \textcolor{black}{substitute} for its gradient in the distributional definition of solution.  With  some abuse of notation, this is the meaning attributed to $\nabla u$ in the statements of Theorems \ref{secondloc}, \ref{secondconvex} and  \ref{seconddir}. 
\\
A  definition of generalized solution to  problem \eqref{eqdirichlet} and to problem \eqref{eqneumann} is given in Section \ref{proofs},  where 
 an existence, uniquess and first-order summability result from \cite{cmapprox} is also recalled. Note that, owing to its uniqueness, this kind of generalized solution agrees with the weak solution whenever $f$ is summable  enough, depending on the nonlinearity of the differential operator, for a weak solution to exist. Generalized local solutions to equation \eqref{localeq} are defined in Section \ref{proofsloc}. 
}
\end{remark}

\section{A differential inequality}\label{diff}

The subject of   this section is  a lower bound for the square of the differential operator on the left-hand side of the equations in \eqref{eqdirichlet}  and \eqref{eqneumann} in terms of  an operator in divergence form, plus (a positive constant times) 
 derivatives of $a(|\nabla u|)\nabla u$ squared. This is a critical step in the proof of our main results, and is the content of the following lemma. 

%

\begin{lemma}\label{lemma1}
Assume that $a \in C^1[0, \infty )$, and that the first inequality
in \eqref{infsup} holds. Then there exists a positive constant $C= C(n, i_a)$ such that 
\begin{align}\label{pointwise}
\big({\rm div} (a(|\nabla u|)\nabla u )\big)^2 & \geq \sum _{j=1}^n
\big(a(|\nabla u|)^2 u_{x_j}\Delta u\big)_{x_j} \\ \nonumber &-
\sum _{i=1}^n \Big(a(|\nabla u|)^2 \sum _{j=1}^n u_{x_j}u_{x_i x_j}\Big)_{x_i}
+ Ca(|\nabla u|)^2 |\nabla ^2 u|^2
\end{align}
for every function $u \in C^3(\Omega )$. Here, $|\nabla ^2u| = \big(\sum _{i,j=1}^n u_{x_ix_j}^2)^{\frac 12}$.
\end{lemma}
\par\noindent
{\bf Proof}.  Let $u \in C^3(\Omega )$. Computations show that
\begin{align}\label{point1}
\big(& {\rm div}  (a(|\nabla u|)\nabla u )\big)^2  = \big(
a(|\nabla u|) \Delta u + a'(|\nabla u|) \nabla |\nabla u| \cdot
\nabla u\big)^2 \\ \nonumber & = a(|\nabla u|)^2 \big((\Delta u)^2
- |\nabla ^2 u|^2\big) + a(|\nabla u|)^2  |\nabla ^2 u|^2 +
 \\ \nonumber & \quad + a'(|\nabla u|)^2 (\nabla |\nabla u| \cdot \nabla u )^2 + 2
a(|\nabla u|) a'(|\nabla u|)\Delta u \nabla |\nabla u| \cdot
\nabla u \\ \nonumber & = a(|\nabla u|)^2 \big(\sum _{j=1}^n
(u_{x_j}\Delta u)_{x_j} - \sum _{i,j=1}^n (u_{x_j} u_{x_ix_j})_{x_i}
\big) + a(|\nabla u|)^2 |\nabla ^2 u|^2
\\ \nonumber & \quad + a'(|\nabla u|)^2 (\nabla |\nabla u| \cdot \nabla u )^2
+ 2 a(|\nabla u|) a'(|\nabla u|)\Delta u \nabla |\nabla u| \cdot
\nabla u
\\ \nonumber & =  \sum _{j=1}^n
(a(|\nabla u|)^2 u_{x_j}\Delta u)_{x_j} - \sum _{i,j=1}^n (a(|\nabla u|)^2 u_{x_j}
u_{x_ix_j})_{x_i}  \\ \nonumber & \quad - 2a(|\nabla u|)
a'(|\nabla u|)\big(\Delta u\, \nabla |\nabla u| \cdot \nabla u -
\sum _{i,j=1}^n |\nabla u|_{x_i}u_{x_j} u_{x_ix_j}\big) \\ \nonumber &
\quad + a(|\nabla u|)^2 |\nabla ^2 u|^2 + a'(|\nabla u|)^2 (\nabla
|\nabla u| \cdot \nabla u )^2 + 2 a(|\nabla u|) a'(|\nabla
u|)\Delta u \nabla |\nabla u| \cdot \nabla u
\\ \nonumber & =
\sum _{j=1}^n
(a(|\nabla u|)^2 u_{x_j}\Delta u)_{x_j} - \sum _{i,j=1}^n (a(|\nabla u|)^2 u_{x_j}
u_{x_ix_j})_{x_i} 
\\ \nonumber & \quad +
2a(|\nabla u|) a'(|\nabla u|) \sum _{i,j=1}^n |\nabla u|_{x_i}u_{x_j}
u_{x_ix_j} + a(|\nabla u|)^2 |\nabla ^2 u|^2 + a'(|\nabla u|)^2
(\nabla |\nabla u| \cdot \nabla u )^2\,,
\end{align}
where $\lq\lq \cdot "$ stands for scalar product in $\rn$. 
After relabeling the indices, one has that
\begin{align}\label{point2}
a'(|\nabla u|)^2 & (\nabla |\nabla u| \cdot \nabla u )^2 +
2a(|\nabla u|) a'(|\nabla u|) \sum _{i,j=1}^n |\nabla u|_{x_i}u_{x_j}
u_{x_ix_j} + a(|\nabla u|)^2 |\nabla ^2 u|^2  \\ \nonumber & =
a(|\nabla u|)^2\bigg[\bigg(\frac {|\nabla u|a'(|\nabla
u|)}{a(|\nabla u|)}\bigg)^2 \bigg(\sum
_{i,k=1}^n\frac{u_{x_k}u_{x_i}}{|\nabla u|^2}u_{x_k x_i}\bigg)^2 \\
\nonumber &  \quad \qquad \quad \quad \quad+ 2\sum _{i,j,k=1}^n\frac {|\nabla u|a'(|\nabla
u|)}{a(|\nabla u|)}\frac{u_{x_k}u_{x_i}}{|\nabla u|^2} u_{x_k
x_j}u_{x_i x_j} + \sum _{i,j=1}^n u_{x_i x_j}^2\bigg].
\end{align}
Now, set $$\omega_u  = \frac{\nabla u}{|\nabla u|},
\quad \vartheta_u = \frac {|\nabla u|a'(|\nabla u|)}{a(|\nabla u|)},\quad
H_u = \nabla ^2 u.$$
Observe that $\omega _u \in \mathbb R^n$, with $|\omega _u|=1$, $H_u$ is a symmettic matrix in $\mathbb R^{n\times n}$, and, by  \eqref{infsup}, $\vartheta _u \geq i_a$.
With this notation in place, the expression in square brackets on the right-hand side of \eqref{point2} takes the form
\begin{align}\label{point3}
\vartheta_u ^2 (H \omega_u \cdot \omega_u )^2 + 2\vartheta_u
H \omega_u \cdot H \omega_u  + {\rm tr }\big(H ^2_u\big)\,,
\end{align}
where $\lq\lq {\rm tr }"$ denotes  the trace of a matrix. The proof of inequality \eqref{pointwise} is thus reduced to showing that 
%
%
\begin{align}\label{point4}
\vartheta_u ^2 (H \omega_u \cdot \omega_u )^2 + 2\vartheta_u
H \omega_u \cdot H \omega_u  + {\rm tr }\big(H ^2_u\big) \geq C {\rm tr }\big(H_u^2\big)
\end{align}
for some   positive constant $C=C(n, i_a)$.
To establish inequality \eqref{point4}, define the function $\psi :  \mathbb R \times \mathbb R^n \times (\mathbb R^{n \times n} \setminus \{0\} )\to \mathbb R$ as
$$\psi (\vartheta, \omega , H) = \vartheta ^2 \frac{(H \omega \cdot \omega )^2}{ {\rm tr } \big(H ^2\big)} + 2\vartheta \frac{H \omega \cdot H
\omega }{ {\rm tr }\big(H ^2\big)} + 1$$
for $(\vartheta , \omega, H) \in \mathbb R \times  \rn \times ( \mathbb R^{n \times n}\setminus \{0\})$, and note that \eqref{point4} will follow if we show that there exists a positive constant $C$ such that
\begin{equation}\label{point4'}
\psi (\vartheta, \omega , H)  \geq C
\end{equation}
if $\vartheta \geq i_a$, $|\omega |=1$ and $H$ is any non-vanishing symmetric matrix $H$.
For each fixed $\omega$ and $H$, the  quadratic function $\vartheta \mapsto \psi (\vartheta, \omega , H)$  attains its 
%
%
%
%
%
 minimum at $\vartheta  = -\, \frac{H \omega \cdot  H
\omega }{(H \omega \cdot \omega )^2}$. We claim that
\begin{equation}\label{point6'}
-\, \frac{H \omega \cdot H \omega }{(H  \omega \cdot \omega )^2}
\leq -1.
\end{equation}
To verify equation \eqref{point6'},  choose a
basis in $\rn$ in which $H$ has diagonal form ${\rm
diag}(\lambda _1, \dots \lambda _n)$, and let  $(\omega _1, \dots , \omega _n)$ denote the vector of the components of $\omega$ with respect to this basis. Then 
$$H \omega \cdot H \omega = \sum _{i=1}^n \lambda _i^2 \omega _i^2\\, \qquad H  \omega \cdot \omega  = \sum _{i=1}^n\lambda _i\omega _i^2\,,$$
whence \eqref{point6'} follows, since
\begin{equation}\label{point6}
\big(\sum _{i=1}^n\lambda _i\omega _i^2 \big)^2\leq \big(\sum _{i=1}^n
\lambda _i^2 \omega _i^2\big) \big(\sum _{i=1}^n\omega _i^2\big) =
\big(\sum _{i=1}^n \lambda _i^2 \omega _i^2\big),
\end{equation}
by Schwarz' inequality. Note that the equality  holds in \eqref{point6} inasmuch as $\sum _{i=1}^n\omega _i^2 = 1$. 
Owing to \eqref{point6'}, $\psi (\vartheta, \omega , H) $ is a stricly increasing function of $\vartheta$ for
$\vartheta \geq -1$. Hence, by the first inequality in \eqref{infsup},
\begin{equation}\label{point100}
\psi (\vartheta, \omega , H)  \geq \psi (i_a, \omega , H) > \psi (-1, \omega , H) 
\end{equation}
if $\vartheta \geq i_a $ and $|\omega |=1$ . Assume, for a moment, that we know that
\begin{equation}\label{point101'}
 \psi (-1, \omega , H)\geq 0
\end{equation}
if $|\omega|=1$ and $H$ is any symmetric matrix. Since $\psi$
 is a continuous function,
we deduce from \eqref{point100} and \eqref{point101'} that
\begin{equation}\label{point102}
\psi (\vartheta, \omega , H)  \geq \psi (i_a, \omega , H) \geq \inf_{|\omega|=1, \,H \,{\rm sym}} \psi (i_a, \omega , H) = \min_{|\omega|=1, \,H  \,{\rm sym}, \, |H|=1} \psi (i_a, \omega , H) >0
\end{equation}
if  $|\omega|=1$ and $H$ is symmetric and different from $0$.  
%
Hence \eqref{point4'} follows. Observe that the equality holds in \eqref{point102} since $\psi$ is a homogenenous function of degree $0$ in $H$.
\\ It remains to prove inequality \eqref{point101'}, namely that
\begin{equation}\label{point104}
(H \omega \cdot \omega )^2 - 2 H \omega \cdot H
\omega  + {\rm tr }\big(H ^2\big) \geq  0
\end{equation}
 if $|\omega|=1$ and $H$ is symmetric.
After diagonalizing $H$ as above, inequality \eqref{point104} reads
\begin{equation}\label{point105}
\sum _{i=1}^n (\omega _i^2 -1)^2 \lambda _i^2 + 2\sum _{1 \leq i <j\leq n} \omega _i^2\omega _j^2 \lambda _i\lambda _j  \geq  0\,,
\end{equation}
if $\sum _{i=1}^n\omega_i^2 =1$ and $\lambda _i \in \mathbb R$ for $i=1, \dots ,n$.  
Inequality \eqref{point105} is a consequence of the following lemma. \qed

\begin{lemma}\label{formaquad}
Assume that $\eta _i \in \mathbb R$ are such that $\eta _i \geq 0$,  $i=1, \dots n$, and $\sum _{i=1}^n \eta_i \leq 1$. Then 
\begin{equation}\label{point101}
\sum _{i=1}^n (\eta_i -1)^2 \lambda _i^2 + 2\sum _{1 \leq i <j\leq n} \eta _i\eta _j \lambda _i\lambda _j  \geq  0
\end{equation}
for every $\lambda _i \in \mathbb R$, $i=1, \dots n$.
\end{lemma}
{\bf Proof}. By Sylvester's criterion, it suffices to show that the determinants of the north-west minors of the $n\times n$ matrix
\begin{equation}\label{quad1}
\begin{pmatrix}
(\eta _1 -1)^2 & \eta _1\eta _2 & \dots  & \eta_1 \eta_n \\
\eta_2\eta_1 & (\eta _2 -1)^2  &  \dots  & \eta _2 \eta _n\\
\vdots  & \vdots & \ddots  &  \vdots  \\
 \eta_n\eta_1  & \eta_n\eta_2   & \dots  & (\eta _n -1)^2
\end{pmatrix}\,,
\end{equation}
associated with the quadratic form on the left-hand side of \eqref{point101}, are nonnegative for every $\eta _i \geq 0$, $i=1, \dots , n$, with $\sum _{i=1}^n \eta_i \leq1$. Since every    minor of this  kind  has the same structure as the entire matrix, and $\sum _{i=1}^2 \eta _i \leq \sum _{i=1}^3 \eta _i  \leq \, \dots \, \leq \sum _{i=1}^n \eta_i \leq 1$, it suffices to prove that just the determinant of the whole matrix in \eqref{quad1} is nonnegative.
To this purpose, let us begin by showing that 
\begin{align}\label{quad2}
\det &\begin{pmatrix}
(\eta _1 -1)^2 & \eta _1\eta _2 & \dots & \eta_1 \eta_n \\
\eta_2\eta_1 & (\eta _2 -1)^2  &  \dots  & \eta _2 \eta _n\\
\vdots  & \vdots & \ddots   & \vdots  \\
 \eta_n\eta_1  & \eta_n\eta_2   & \dots   & (\eta _n -1)^2
\end{pmatrix}
\\ & = \eta_1^2 (1-2\eta _2)(1-2\eta_3) \times \cdots \times (1-2\eta_n) + \eta_2^2 (1-2\eta _1)(1-2\eta_3) \times \cdots \times (1-2\eta_n) +  \cdots \nonumber
\\ & \quad \cdots +\eta_n^2 (1-2\eta _1)(1-2\eta_2) \times \cdots \times (1-2\eta_{n-1}) + (1-2\eta _1)(1-2\eta_2) \times \cdots \times (1-2\eta_n). \nonumber
\end{align}
Equation \eqref{quad2} can be verified by induction on $n$. The case when $n=2$ is trivial. Assume that  \eqref{quad2} holds with $n$ replaced by $n-1$. We have that
\begin{multline}\label{quad3}
\det \begin{pmatrix}
(\eta _1 -1)^2 & \eta _1\eta _2 & \dots &  \eta_1 \eta_n \\
\eta_2\eta_1 & (\eta _2 -1)^2  &  \dots  & \eta _2 \eta _n\\
\vdots  & \vdots & \ddots  & \vdots  \\
 \eta_n\eta_1  & \eta_n\eta_2   & \dots   & (\eta _n -1)^2
\end{pmatrix}
\\
= \det \begin{pmatrix}
\eta _1 ^2 & \eta _1\eta _2 & \dots & \eta_1 \eta_n \\
\eta_2\eta_1 & (\eta _2 -1)^2  &  \dots & \eta _2 \eta _n\\
\vdots  & \vdots & \ddots  & \vdots   \\
 \eta_n\eta_1  & \eta_n\eta_2 & \dots  & (\eta _n -1)^2
\end{pmatrix} + \,
\det \begin{pmatrix}
(1-2\eta _1 ) & \eta _1\eta _2 & \dots & \eta_1 \eta_n \\
0 & (\eta _2 -1)^2  &  \dots  & \eta _2 \eta _n\\
\vdots  & \vdots & \ddots   & \vdots  \\
 0  & \eta_n\eta_2   & \dots   & (\eta _n -1)^2
\end{pmatrix}
\end{multline}
Our induction assumption tells us that
\begin{align}\label{quad4}
\det &\begin{pmatrix}
(1-2\eta _1 ) & \eta _1\eta _2 & \dots & \eta_1 \eta_n \\
0 & (\eta _2 -1)^2  &  \dots & \eta _2 \eta _n\\
\vdots  & \vdots & \ddots  & \vdots  \\
 0  & \eta_n\eta_2   & \dots   & (\eta _n -1)^2
\end{pmatrix}
\\ \nonumber & = (1-2\eta _1)  \det \begin{pmatrix}
 (\eta _2 -1)^2  &  \dots  & \eta _2 \eta _n\\
   \vdots & \ddots   & \vdots  \\
  \eta_n\eta_2   & \dots &  (\eta _n -1)^2
\end{pmatrix}
\\
\nonumber
& 
= \eta_2^2 (1-2\eta _1)(1-2\eta_3) \times \cdots \times (1-2\eta_n) + \cdots  + \eta_n^2 (1-2\eta _1)(1-2\eta_2)  \times \cdots \times (1-2\eta_{n-1})  
\\ \nonumber & \quad 
+ (1-2\eta _1)(1-2\eta_2) \times \cdots \times (1-2\eta_n). 
\end{align}
On the other hand, we claim that
\begin{equation}\label{quad6}
\det  \begin{pmatrix}
\eta _1 ^2 & \eta _1\eta _2 & \dots &  \eta_1 \eta_n \\
\eta_2\eta_1 & (\eta _2 -1)^2  &  \dots & \eta _2 \eta _n\\
\vdots  & \vdots & \ddots  & \vdots  \\
 \eta_n\eta_1  & \eta_n\eta_2   & \dots  & (\eta _n -1)^2
\end{pmatrix}
 = \eta_1^2 (1-2\eta_2)(1-2\eta _3) \cdots (1-2\eta _n).
\end{equation}
Equation \eqref{quad6} can be proved by induction again. If $n=2$, this equation can be verified via a direct computation. Assume now that it holds with $n$ replaced by $(n-1)$. Then, 
\begin{align}\label{quad5}
\det & \begin{pmatrix}
\eta _1 ^2 & \eta _1\eta _2 &  \eta _1\eta _3 & \dots & \eta_1 \eta_n \\
\eta_2\eta_1 &  (\eta _2 -1)^2   &  \eta _2\eta _3 & \dots& \eta _2 \eta _n\\
 \eta_3\eta_1  & \eta_3\eta_2  & (\eta _3 -1)^2 & \dots   &  \eta _3\eta _n\\
\vdots & \vdots & \vdots & \ddots & \vdots \\
 \eta_n\eta_1  & \eta_n\eta_2   &  \eta _n\eta _3 & \dots &  (\eta _n -1)^2
\end{pmatrix}
%
\\ \nonumber
& 
= 
\det  \begin{pmatrix}
\eta _1 ^2 & \eta _1\eta _2 &  \eta _1\eta _3 & \dots & \eta_1 \eta_n \\
\eta_2\eta_1 & \eta _2^2  &  \eta _2\eta _3 & \dots& \eta _2 \eta _n\\
\eta_3\eta_1  & \eta_3\eta_2  & (\eta _3 -1)^2 & \dots   &  \eta _3\eta _n\\
\vdots & \vdots & \vdots & \ddots & \vdots \\
 \eta_n\eta_1  & \eta_n\eta_2   &  \eta _n\eta _3 & \dots &  (\eta _n -1)^2
\end{pmatrix} 
\\ \nonumber & \quad 
+
\det  \begin{pmatrix}
\eta _1 ^2 & 0 &  \eta _1\eta _3 & \dots   & \eta_1 \eta_n \\
\eta_2\eta_1 & 1- 2 \eta _2  &  \eta _2\eta _3 &  \dots &\eta _2 \eta _n\\
\eta_3\eta_1  & 0  & (\eta _3 -1)^2 & \dots   &  \eta _3\eta _n\\
\vdots & \vdots & \vdots & \ddots & \vdots \\
 \eta_n\eta_1  & 0   &  \eta _n\eta _3 & \dots &  (\eta _n -1)^2
\end{pmatrix}
\\ \nonumber
& 
= 0 +  (1-2\eta _2)\det  \begin{pmatrix}
\eta _1 ^2 & \eta _1\eta _3 & \dots &  \eta_1 \eta_n \\
\eta_3\eta_1  &  (\eta _3 -1)^2 & \dots  &  \eta _3\eta _n\\
\vdots & \vdots & \ddots & \vdots \\
 \eta_n\eta_1  & \eta_n\eta_3   &   \dots &  (\eta _n -1)^2
\end{pmatrix}
\\ \nonumber
& 
=   \eta_1^2 (1-2\eta_2) (1-2\eta_3) \cdots (1-2\eta_n).
\end{align}
Note that in the last equality we have made use of the  induction assumption, and of the  fact that the determinant of a matrix with a couple of linearly dependent columns vanishes.
%
%
Equation \eqref{quad2} follows from \eqref{quad3}, \eqref{quad4} and \eqref{quad6}.
\\ With equation \eqref{quad2} at disposal, let us define the function $\phi : \mathbb R^n \to \mathbb R$ as 
\begin{align}\label{quad8}
\phi(\eta) & =
\eta_1^2 (1-2\eta_2) \times  \cdots \times (1-2\eta _n) + \cdots + \eta_n^2 (1-2\eta_1) \times \cdots \times (1-2\eta _{n-1}) \\ \nonumber & \quad + (1-2\eta_1) (1-2\eta _2) \times \cdots  \times (1-2\eta _{n}) 
\end{align}
for $\eta \in \rn$, where we have set $\eta =(\eta _1, \dots , \eta _n)$.
Define
$$A= \Big\{\eta \in \rn: \eta _i \geq 0, \, i=1, \dots , n, \, \sum_{i=1}^n\eta_i \leq 1 \Big\}.$$
We have to show that
\begin{equation}\label{quad15}
\phi (\eta) \geq 0\qquad \hbox{for every $\eta \in A$ .}
\end{equation}
On performing the products on the right-hand side of \eqref{quad8}, and rearranging the resulting terms, one can verify that

\begin{align}\label{quad11}
\phi(\eta) & =  \eta _1^2 \bigg[1 + (-2) \sum_{i\neq 1} \eta_i + (-2)^2 \sum _
{\substack{
\scriptstyle i_1<i_2\\
\scriptstyle i_1,i_2\neq 1} }\eta_{i_1}\eta_{i_2} + \cdots 
\\ \nonumber & \quad \quad \quad 
\cdots +   (-2)^k \sum _
{\substack{
\scriptstyle i_1<i_2<\cdots<i_k\\
\scriptstyle i_1, \cdots , i_k \neq 1} }
 \eta_{i_1}\eta_{i_2} \cdots \eta_{i_k} + \cdots 
+ (-2)^{n-1}\eta_{2} \cdots \eta_{n}\bigg] + \cdots
\\ \nonumber & \qquad \qquad \qquad  \qquad \qquad \qquad \vdots
\\ \nonumber & \quad  \cdots+  \eta _n^2 \bigg[1 + (-2) \sum_{i\neq n} \eta_i + (-2)^2 \sum _
{\substack{
\scriptstyle i_1<i_2\\
\scriptstyle i_1,  i_2 \neq n} } \eta_{i_1}\eta_{i_2} + \cdots
\\ \nonumber & \quad  \quad \quad 
 \cdots +    (-2)^k\sum _{\substack{
\scriptstyle i_1<i_2<\cdots<i_k\\
\scriptstyle i_1, \cdots , i_k \neq n} }\eta_{i_1}\eta_{i_2} \cdots \eta_{i_k} + \cdots 
+ (-2)^{n-1}\eta_{1} \cdots \eta_{n-1}\bigg]
\\ \nonumber & \quad +
1 + (-2)\sum _{i=1,\cdots , n} \eta_{i} + (-2)^2\sum _{i_1<i_2} \eta_{i_1}\eta_{i_2}  + (-2)^3  \sum _{i_1<i_2<i_3} \eta_{i_1}\eta_{i_2} \eta_{i_3} + \cdots  + (-2)^n \eta_1 \cdots \eta_n \,.
\end{align}
Let us denote by $\mathcal S_k$, for $k= 1, \dots , n$, the elementary symmetric functions of the $n$ numbers $\eta_1, \dots , \eta_n$. Namely,
$$\mathcal S_k = \sum _{i_1<i_2<\cdots<i_k} \eta_{i_1}\eta_{i_2} \cdots \eta_{i_k}.$$
Observe that
\begin{align}\label{feb30}
(1-\mathcal S_1)^2  = \Big(1 - \sum_{i=1}^n \eta _i\Big)^2 = 1 - 2\sum_{i=1}^n \eta _i
+ 2 \sum _{i_1<i_2} \eta_{i_1}\eta_{i_2} + \sum_{i=1}^n \eta _i^2\,.
\end{align}
Moreover,
\begin{align}\label{feb31}
\mathcal S_1 \mathcal S_k & = \sum_{i=1, \dots n} \eta _i   \sum _{i_1<\dots<i_k} \eta_{i_1}\eta_{i_2} \cdots \eta_{i_k}
\\ \nonumber & = \eta _1 ^2 \sum _
{\substack{
\scriptstyle i_1<i_2<\cdots<i_{k-1}\\
\scriptstyle i_1, \cdots , i_{k-1} \neq 1} }
 \eta_{i_1}\eta_{i_2} \cdots \eta_{i_{k-1}} +\,  \cdots  
\, + \eta _n ^2 \sum _
{\substack{
\scriptstyle i_1<i_2<\cdots<i_{k-1}\\
\scriptstyle i_1, \cdots , i_{k-1} \neq n} } 
 \eta_{i_1}\eta_{i_2} \cdots \eta_{i_{k-1}} 
\\ \nonumber & \qquad
+ (k+1) \sum _
{ i_1<i_2<\cdots<i_{k+1} } 
 \eta_{i_1}\eta_{i_2} \cdots \eta_{i_{k+1}}
\\ \nonumber &= \eta _1 ^2 \sum _
{\substack{
\scriptstyle i_1<i_2<\cdots<i_{k-1}\\
\scriptstyle i_1, \cdots , i_{k-1} \neq 1} }
 \eta_{i_1}\eta_{i_2} \cdots \eta_{i_{k-1}} +\,  \dots  \, + \eta _n ^2 \sum _
{\substack{
\scriptstyle i_1<i_2<\cdots<i_{k-1}\\
\scriptstyle i_1, \cdots , i_{k-1} \neq n} } 
 \eta_{i_1}\eta_{i_2} \cdots \eta_{i_{k-1}} + (k+1) \mathcal S_{k+1}
\end{align}
for $k=2, \dots , n-1$, and 
\begin{align}\label{feb38}
\mathcal S_1 \mathcal S_n  = \bigg(\sum_{i=1, \dots n} \eta _i \bigg)  \eta_1 \,\cdots \, \eta_n
= \eta _1 ^2 \eta _2 \, \cdots \, \eta _n +\,  \dots  \, + \eta _n ^2 \eta_1 \, \cdots \, \eta _{n-1}\,.
\end{align}
\color{black}
 On making use of equations \eqref{feb30}, \eqref{feb31} and \eqref{feb38}, one can combine the terms on the right-hand side of equation \eqref{quad11} and infer that
\begin{align}\label{feb32}
\phi (\eta) & = (1-\mathcal S _1) \bigg[1 + \sum _{k=1}^{\textcolor{black}{n}} (-1)^k2^{k-1}\mathcal S_k\bigg] + \sum _{k=3}^n (-1)^{k-1}(k-2)2^{k-2}\mathcal S_k.
\end{align}
Since $\mathcal S_1 = \sum_{i=1}^n \eta_i$,  we have that
\begin{equation}\label{feb10}
1 - \mathcal S_1 \geq 0 \quad \hbox{for $\eta \in A$.}
\end{equation}
The sums on the right-hand  side of  equation \eqref{feb32} can be estimated from below via the inequality 
\begin{equation}\label{quad14}
\mathcal S_{k+1} \leq \frac {n-k}{n(k+1)}\mathcal S_k \mathcal S_1 \leq   \frac {n-k}{n(k+1)}\mathcal S_k \quad \hbox{for $\eta \in A$,}
\end{equation}
 and $k=1, \dots , n-1$.
Note that the second inequality in \eqref{quad14} holds by \eqref{feb10}, whereas the first one follows via an iterated use of Newton's inequality \cite[Theorem 51]{HLP}.
%
%
We claim that 
\begin{equation}\label{feb11}
1 + \sum _{k=1}^{n} (-1)^k2^{k-1}\mathcal S_k = 1 - \mathcal S_1 + \sum _{k=2}^{n} (-1)^k2^{k-1}\mathcal S_k \geq 0 \quad  \hbox{for $\eta \in A$.}
\end{equation}
Indeed, by \eqref{quad14},
\begin{equation}\label{feb12}
2^{2h-1}\mathcal S_{2h} - 2^{2h}\mathcal S_{2h+1} \geq 0,
\end{equation}
if $1 \leq h \leq \tfrac {n-1}2$. When $n$ is odd, the
sum starting from $k=2$ in \eqref{feb11} is exhausted by differences of the form appearing in \eqref{feb12}.  When $n$ is even,  this sum contains an additional nonnegative term. Hence, inequality \eqref{feb11} follows.
We next observe that
%
\begin{align}\label{quad12}
 \sum _{k=3}^n (-1)^{k-1}(k-2)2^{k-2}\mathcal S_k \geq 0 \quad  \hbox{for $\eta \in A$.}
\end{align}
%
Actually, inequality \eqref{quad14}  again ensures that
\begin{align}\label{quad50}
(2h-1)2^{2h-1} \mathcal S_{2h+1} -   2h2^{2h} \mathcal S_{2h+2}\geq 0 ,
\end{align}
if $1 \leq h \leq   \tfrac {n-2}2$. When $n$ is even, the sum in \eqref{quad12}  is exhausted by differences of the form appearing in \eqref{quad50}. 
%
When $n$ is odd, this sum contains an additional nonnegative term. Inequality \eqref{quad12} is thus established. Inequality \eqref{quad15} follows from \eqref{feb32}, via \eqref{feb10}, \eqref{feb11} and \eqref{quad12}. 
 Note that, in fact,  
$$
\min_{\eta \in A}\phi (\eta) =0,
$$
inasmuch as $\phi (\eta)=0$ whenever $\eta$ is a vector all of whose components vanish, but just one, and the latter equals one.
%
%
The proof is complete.  
\qed

\section{Global estimates}\label{proofs}
This section is devoted to proving  Theorems \ref{secondconvex} and \ref{seconddir}. As a preliminary, we briefly discuss the notion of generalized solutions adopted in our results, and recall some of their basic  properties.
\par When the function $f$ appearing on the right-hand side of the equation in problems \eqref{eqdirichlet} or \eqref{eqneumann} has a sufficiently high degree of summability to belong to the dual of the Sobolev type space associated with the function $a$,  weak solutions to the relevant problems are well defined.   In particular, the existence and uniqueness of these solutions can be established via standard monotonicity methods. We are not going to give details in this connection, since they are not needed for our purposes, and  refer the interested reader to \cite{cmapprox} for an account on this issue. We rather focus on the case when $f$ merely belongs  to $L^q(\Omega)$ for any $q\geq 1$. A definition of generalized solution in this case 
involves the use of spaces that consist of  functions whose truncations are weakly differentiable. Specifically, given any $t>0$,  let  $T_{t} : \R \rightarrow \R$ denote the function
defined as $
T_{t} (s) = 
s$  if $|s|\leq t$,  and  $T_{t} (s)=
 t \,{\rm sign}(s) $ if $|s|>t$.
%
%
%
We set
\begin{equation}\label{503}
\mathcal T^{1,1}_{\rm loc}(\o) = \left\{u \, \hbox{is measurable in $\o$}: \hbox{$T_{t}(u)
\in W^{1,1}_{\rm loc}(\o)$ for every $t >0$} \right\}.
\end{equation}
The spaces  $\mathcal T^{1,1}(\o)$ and    $\mathcal T^{1,1}_0(\o)$ are defined accordingly, on
replacing $W^{1,1}_{\rm loc}(\o)$ with $W^{1,1}(\o)$  and $W_{0}^{1,1}(\o)$, respectively, on the right-hand
side of \eqref{503}.
\\ If $u \in\mathcal T^{1,1}_{\rm loc}(\o)$, there exists
a (unique) measurable function  $Z_u : \o \to \rn$ such that
\begin{equation}\label{504}
\nabla \big(T_{t}(u)\big) = \chi _{\{|u|<t\}}Z_u  \qquad \quad
\hbox{a.e. in $\o$}
\end{equation}
for every $t>0$ -- see \cite[Lemma 2.1]{BBGGPV}. Here $\chi _E$ denotes
the characteristic function of the set $E$. 
As already mentioned in Section \ref{sec1}, with abuse of
notation, for every $u \in \mathcal T^{1,1}_{\rm loc}(\o)$  we denote $Z_u$ simply
by $\nabla u$.

%

\smallskip
\par
Assume  that $f \in L^q(\o)$ for some $q \geq 1$.  A function \textcolor{black}{$u \in \mathcal T ^{1,1}_0(\o)$} will be called a generalized solution to the Dirichlet problem \eqref{eqdirichlet} if $a(|\nabla u|)\nabla u \in L^1(\o)$,
\begin{equation}\label{231}
\int_\o a(|\nabla u|) \nabla u \cdot \nabla \varphi \, dx=\int_\o f \varphi \,dx\,
\end{equation}
for every $\varphi \in C^{\infty}_0(\o)$, and there exists a sequence 
 $\{f_k\}\subset  C^\infty _0(\o)$
such that 
$f_k \to f$   in $L^q(\o)$ and the sequence of weak solutions  $\{u_k\}$   to the problems \eqref{eqdirichlet} with $f$ replaced by $f_k$ satisfies 
%
%
$$u_k\to u \quad
\hbox{ a.e. in $\o$.}$$
In \eqref{231}, $\nabla u$ stands for the function $Z_u$ fulfilling \eqref{504}.
\\
By \cite{cmapprox},  there exists a unique generalized solution $u$ to problem \eqref{eqdirichlet},  and
\begin{equation}\label{estgraddir}
\|a(|\nabla u|) \nabla u\|_{L^1(\o)} \leq C \|f\|_{L^1(\o)} 
\end{equation} 
for some constant  $C=C(|\o|, n, i_a, s_a)$.  Moreover, if $\{f_k\}$ is any sequence as above, and $\{u_k\}$ is the associated sequence of weak solutions, then 
\begin{equation}\label{gen1}
u_k \to u \quad \hbox{and} \quad \nabla u_k \to \nabla u \quad \hbox{a.e. in $\Omega$,}
\end{equation}
up to subsequences.

\smallskip
\par 
The definition of generalized solutions to the Neumann problem \eqref{eqneumann} can be given analogously. 
Assume that
 $f \in L^q(\o)$ for some $q \geq 1$, and satisfies \eqref{0mean}. 
A function \textcolor{black}{$u \in \mathcal T ^{1,1}(\o)$} will be called a  generalized solution to problem \eqref{eqneumann} if $a(|\nabla u|)\nabla u\in L^1(\o)$, equation \eqref{231} holds for every   $\varphi \in C^\infty (\o) \cap W^{1,\infty}(\o)$, and
there exists a sequence $\{f_k\}\subset  C^\infty _0(\o)$,
with  $\int _\o f_k(x)\, dx =0$ for $k \in \N$, such that
$f_k \to f$   in $L^q(\o)$ and the sequence of  (suitably normalized by additive constants) weak solutions  $\{u_k\}$   to the problems \eqref{eqneumann} with $f$ replaced by $f_k$ satisfies 
%
%
$$u_k\to u \quad
\hbox{ a.e. in $\o$.}$$
%
Owing to \cite{cmapprox}, if $\Omega$ is a bounded Lipschitz domain, then
there exists a unique (up to addive constants) generalized solution $u$ to problem \eqref{eqneumann},  and
\begin{equation}\label{estgradneu}
\|a(|\nabla u|) \nabla u\|_{L^1(\o)} \leq C \|f\|_{L^1(\o)} 
\end{equation} 
for some constant  $C=C( L_\Omega, d_\o, n,  i_a, s_a)$.    Moreover,
$\{f_k\}$ is any sequence as above, and $\{u_k\}$ is the associated sequence of (normalized) weak solutions, then 
\begin{equation}\label{gen1neu}
u_k \to u \quad \hbox{and} \quad \nabla u_k \to \nabla u \quad \hbox{a.e. in $\Omega$,}
\end{equation}
up to subsequences.
\par

We conclude our background by recalling the definitions of Marcinkiewicz, and, more generally,  Lorentz spaces that enter in our results.  Let $(\MR, m)$ be a $\sigma$-finite non atomic measure space. Given $q \in [1,\infty]$, the Marcinkiewicz space $L^{q, \infty} (\MR, m)$, also called weak $L^q(\MR, m)$ space, is the Banach function space endowed with the norm defined as 
\begin{equation}\label{weakleb}
\|\psi\|_{L^{q, \infty} (\MR, m)} = \sup _{s \in (0, m(\MR))} s ^{\frac 1q} \psi^{**}(s)
\end{equation}
for a measurable function $\psi$ on $\MR$. Here, $\psi^*$ denotes the decreasing rearrangement of $\psi$, and $\psi ^{**}(s)= \smallint _0^s \psi^* (r)\, dr$ for $s >0$.
The space $L^{q, \infty} (\MR, m)$ is borderline in the family of Lorentz spaces $L^{q, \sigma}(\MR )$, with  $q \in [1, \infty ]$ and $\sigma \in [1, \infty ]$, that are equipped with the norm given by
\begin{equation}\label{lorentz}
\|\psi\|_{L^{q, \sigma}(\MR )} = \|s^{\frac 1q - \frac 1\sigma } \psi^{**}(s)
\|_{L^\sigma (0, m (\MR))}
\end{equation}
for $\psi$ as above.
Indeed, one has that 
\begin{equation}\label{lorentz11}
  L^{q, \sigma _1}(\MR )
\subsetneqq L^{q, \sigma _2}(\MR ) \qquad \hbox{if $q \in [1, \infty ]$ and
$1 \leq \sigma _1 < \sigma _2\leq \infty$}.
\end{equation}
Also
$$L^{q,q}(\MR) = L^q(\MR) \qquad \hbox{for $q \in (1,\infty]$,}$$
up to equivalent norms.
In the limiting case when $q=1$, the Marcinkiewicz type space $L^{1, \infty} \log L (\MR, m)$ comes into play in our results as a replacement for $L^{1, \infty} (\MR, m)$, which agrees with $L^1(\MR, m) $. A norm in $L^{1, \infty} \log L (\MR, m)$  is defined as 
\begin{equation}\label{weaklog}
\|\psi\|_{L^{1, \infty} \log L (\MR, m)} = \sup _{s \in (0, m(\MR))} s \log\big(1+ \tfrac{C}s\big) \psi^{**}(s),
\end{equation}
for any constant  $C>m(\MR)$.  Different constants $C$ result in equivalent norms in \eqref{weaklog}.

\smallskip
\noindent
{\bf Proof of Theorem \ref{seconddir}}. We begin with a proof in the case when $u$ is the generalized solution to the Dirichlet problem \eqref{eqdirichlet}. The needed variants for the solution to the Neumann problem \eqref{eqneumann} are indicated at the end.
\\
The proof is split in steps. In Step 1 we establish  the result under some additional regularity assumptions on $a$, $\Omega$ and $f$.   The remaining steps are devoted to removing the extra  assumptions, by approximation.
\\ \emph{Step 1}. Here, we assume that the  following extra conditions are in force:
\begin{equation}\label{fsmooth}
f \in C^\infty _0 (\Omega);
\end{equation}
 \begin{equation}\label{omegasmooth}
\partial \Omega \in C^\infty;
\end{equation}
\color{black}
\begin{equation}\label{abound}
 a: [0, \infty) \to [0, \infty) \quad \hbox{and} \quad c_1 \leq a(t) \leq c_2 \quad \hbox{for $t \geq 0$, }
\end{equation}
for some constants $c_2 > c_1 >0$; the function $\mathcal A : \rn \to [0, \infty)$, defined as $\mathcal A (\eta) = a(|\eta|)$ for $\eta \in \rn$, is such that 
\begin{equation}\label{Asmooth}
\mathcal A  \in C^\infty(\rn).
\end{equation}
Standard regularity results then ensure that the solution $u$ to problem \eqref{eqdirichlet} is classical, and $u \in C^\infty(\overline \o)$ (see e.g. \cite[Proof of Theorem 1.1]{cmCPDE} for details). 
\color{black} Let $\xi \in C^\infty _0(\rn)$. Squaring both sides of the equation in \eqref{eqdirichlet}, multiplying through the resulting equation by $\xi^2$, integrating both sides over $\Omega$, and making use of   inequality  \eqref{pointwise} yield
\begin{align}\label{main1}
\int _\Omega \xi^2 f^2\, dx &= \int _\Omega \xi^2 \big({\rm div}(a(|\nabla u|)\nabla u)\big)^2\, dx 
\\ \nonumber & \geq \int_\Omega \xi ^2 \bigg[\sum _{j=1}^n
\big(a(|\nabla u|)^2 u_{x_j}\Delta u\big)_{x_j} -
\sum _{i=1}^n \Big(a(|\nabla u|)^2 \sum _{j=1}^n u_{x_j}u_{x_i x_j}\Big)_{x_i}\bigg]\, dx 
\\ \nonumber & \quad + 
C \int _\Omega \xi ^2 a(|\nabla u|)^2 |\nabla ^2 u|^2\, dx\,
\end{align}
for some constant $C =C (n, i_a)$.
Now, \cite[Equation
(3,1,1,2)]{Grisvard} tells us that
\begin{multline}\label{grisv1}
\Delta u \frac{\partial u}{\partial \nu } - \sum _{i,j=1}^n
u_{x_i x_j} u_{x_i} \nu _j \\ = {\rm div }_T
\bigg(\frac{\partial u}{\partial \nu } \nabla _T u\bigg) -
{\rm tr}\mathcal B \bigg(\frac{\partial u}{\partial \nu }\bigg)^2
- \mathcal B (\nabla _T \, u , \nabla _T \, u) - 2 \nabla _T\,
u \cdot \nabla _T\, \frac{\partial u}{\partial \nu }
 \qquad \quad
\hbox{on $\partial \Omega$,}
\end{multline}
where $\mathcal B$ is the second fundamental form on $\partial
\Omega$, ${\rm tr}\mathcal B $ is its trace,  
 ${\rm div }_T$ and
$\nabla _T $ denote the divergence and the gradient
operator on $\partial \Omega$, respectively, and $\nu _j$ stands for the $j$-th component of $\nu$. From 
 the divergence theorem and equation \eqref{grisv1} we deduce that 
\begin{align}\label{main2}
\int_\Omega \xi ^2&  \Big[\sum _{j=1}^n
\big(a(|\nabla u|)^2 u_{x_j}\Delta u\big)_{x_j} -
\sum _{i=1}^n \big(a(|\nabla u|)^2 \sum _{j=1}^n u_{x_j}u_{x_i x_j}\big)_{x_i}\Big]\, dx 
\\ \nonumber & = \int _{\partial \Omega } \xi^2 a(|\nabla u|)^2 \Big[\Delta u \frac {\partial u}{\partial \nu} - \sum _{i,j=1}^nu_{x_i x_j}u_{x_i}\nu _j \Big]\, d\mathcal H^{n-1}(x) 
\\ \nonumber & \quad  \quad- 2 \int _\Omega  a(|\nabla u|)^2\xi \nabla \xi  \cdot \Big[\Delta u \nabla u - \sum _{j=1}^n u_{x_j}\nabla u_{x_j}\Big]\, dx
\\ \nonumber & = \int _{\partial \Omega } \xi^2 a(|\nabla u|)^2 \bigg[{\rm div }_T
\bigg(\frac{\partial u}{\partial \nu } \nabla _T u\bigg) -
{\rm tr}\mathcal B \bigg(\frac{\partial u}{\partial \nu }\bigg)^2
\\ \nonumber & \quad \quad  \qquad \qquad \qquad 
- \mathcal B (\nabla _T \, u , \nabla _T \, u) - 2 \nabla _T\,
u \cdot \nabla _T\, \frac{\partial u}{\partial \nu }\bigg] \,d\mathcal H^{n-1}(x)
\\ \nonumber & \quad - 2 \int _\Omega  a(|\nabla u|)^2\xi \nabla \xi  \cdot \Big[\Delta u \nabla u - \sum _{j=1}^n u_{x_j}\nabla u_{x_j}\Big]\, dx\,.
\end{align}
By Young's inequality, there exists a constant $C=C(n)$ such that
\begin{align}\label{main4}
2\bigg|\int _\Omega  & a(|\nabla u|)^2\xi \nabla \xi  \cdot \Big[\Delta u \nabla u - \sum _{j=1}^n u_{x_j}\nabla u_{x_j}\Big]\, dx\bigg| 
\\ \nonumber  & \leq \varepsilon C \int _\Omega \xi^2 a(|\nabla u|)^2 |\nabla^2 u|^2\,dx + \frac C \varepsilon \int _\Omega |\nabla \xi|^2 a(|\nabla u|)^2 |\nabla u|^2\,dx\,
\end{align}
for every  $\varepsilon >0$.
Equations \eqref{main1}, \eqref{main2} and \eqref{main4} ensure that there exist constants $C=C(n, i_a)$ and $C'=C'(n, i_a)$ such that
\begin{align}\label{main3bis}
C(1-\varepsilon) \int _\Omega \xi^2 a(|\nabla u|)^2 |\nabla^2 u|^2\,dx & \leq 
\int _\Omega \xi^2 f^2\, dx  + \frac {C'} \varepsilon \int _\Omega |\nabla \xi|^2 a(|\nabla u|)^2 |\nabla u|^2\,dx 
\\ \nonumber & \quad + 
\bigg|\int _{\partial \Omega } \xi^2 a(|\nabla u|)^2 \bigg[{\rm div }_T
\bigg(\frac{\partial u}{\partial \nu } \nabla _T u\bigg) -
{\rm tr}\mathcal B \bigg(\frac{\partial u}{\partial \nu }\bigg)^2
\\ \nonumber & \quad \quad \quad
- \mathcal B (\nabla _T \, u , \nabla _T \, u) - 2 \nabla _T\,
u \cdot \nabla _T\, \frac{\partial u}{\partial \nu }\bigg] \,d\mathcal H^{n-1}(x)\bigg|\,.
\end{align}
On the other hand, owing to the Dirichlet boundary condition,  $\nabla _Tu=0$ on $\partial \Omega$, and hence  
\begin{align}\label{main3}
\bigg|\int _{\partial \Omega }& \xi^2 a(|\nabla u|)^2 \bigg[{\rm div }_T
\bigg(\frac{\partial u}{\partial \nu } \nabla _T u\bigg) -
{\rm tr}\mathcal B \bigg(\frac{\partial u}{\partial \nu }\bigg)^2
\\ \nonumber & \quad \quad \quad \quad \quad \quad
- \mathcal B (\nabla _T \, u , \nabla _T \, u) - 2 \nabla _T\,
u \cdot \nabla _T\, \frac{\partial u}{\partial \nu }\bigg] \,d\mathcal H^{n-1}(x)\bigg| 
\\ \nonumber & = \bigg|- 
\int _{\partial \Omega } \xi^2 a(|\nabla u|)^2 
{\rm tr}\mathcal B \bigg(\frac{\partial u}{\partial \nu }\bigg)^2 \,d\mathcal H^{n-1}(x)\bigg| 
\\ \nonumber & \leq C \int _{\partial \Omega }\xi ^2a(|\nabla u|)^2 |\nabla u|^2 |\mathcal B| \,d\mathcal H^{n-1}(x)\,,
\end{align}
for some constant $C=C(n)$. Here, $|\mathcal B|$ denotes the norm of $\mathcal B$. Next, assume that
\begin{equation}\label{feb22}
\xi \in C^\infty _0(B_r(x))
\end{equation}
for some $x \in \overline \Omega$ and $r>0$. 
\\ First, suppose that
$x \in \partial \o$. 
Let us distinguish the cases when $n \geq 3$ or $n=2$. When $n \geq 3$,
set
\begin{equation}\label{main5}
Q(r) = \sup _{x \in \partial \o} \sup _{E \subset \partial \Omega \cap B_r(x)} \frac{\int _E |\mathcal B| \,d \hh (y)}{{\rm cap}(E)}\, \quad\quad \hbox{for $r>0$,}
\end{equation}
where ${\rm cap}(E)$ stands for the  capacity of the set
 $E$  given by
\begin{equation}\label{cap}
{\rm cap}(E) = \inf \bigg\{\int _\rn |\nabla v|^2\, dy : v\in C^1_0(\rn), v\geq 1 \,\hbox{on }\, E\bigg\}.
\end{equation}
\textcolor{black}{A weighted trace inequality on half-balls \cite{Mazya62,  Mazya64} (see also \cite[Section 2.5.2]{Mabook})},  combined with a local flattening argument for $\Omega$ on a half-space, and with an even-extension argument from a half-space into $\rn$, ensures that there exists a constant $C=C(L_\Omega,  d_\o, n)$ such that
\begin{equation}\label{main6}
\int _{\partial \Omega \cap B_r(x)} v^2 \, |\mathcal B| \, d \hh (y) \leq C \textcolor{black}{Q(r)}
\int _{\Omega \cap B_r(x)}|\nabla v|^2 \, dy
\end{equation}
for every $x \in \partial \Omega$, $r >0$ and $v \in C^1_0(B_r(x))$. 
Furthermore,
a standard trace inequality  tells us that that there exists a constant  $C=C(L_\Omega,  d_\o, n)$  such that
\begin{equation}\label{traceineq}
\bigg(\int _{\partial \Omega \cap B_r(x)} |v|^{\frac{2(n-1)}{n-2}} \, d \hh (y) \bigg)^{\frac{n-2}{n-1}} \leq C \int _{\Omega \cap B_r(x)}|\nabla v|^2 \, dy
\end{equation}
for every  $x \in \partial \Omega$, $r >0$ and $v \in C^1_0(B_r(x))$. By definition \eqref{cap}, choosing trial functions $v$ in \eqref{traceineq} such that $v\geq 1$ on $E$ implies that
\begin{equation}\label{main8}
\hh (E)^{\frac {n-2}{n-1}} \leq C \, {\rm cap} (E) 
\end{equation}
for every set $E \subset \partial \Omega$. 
By a basic property of the decreasing rearrangement (with respect to $\mathcal H^{n-1}$)  \cite[Chapter 2, Lemma 2.1]{BS}, and  \eqref{main8},
\begin{align}\label{main9}
Q(r) & \leq  \sup _{x \in \partial \o} \sup _{E\subset \partial \o \cap B_r(x)} \frac{\int _0^{\hh (E)} (|\mathcal B| _{|\partial \Omega \cap B_r(x)})^*(r)\,dr}{{\rm cap}(E)}
\\ \nonumber &\leq C \sup _{x \in \partial \o} \sup _{s>0} \frac{\int _0^{s} (|\mathcal B| _{|\partial \Omega \cap B_r(x)})^*(r)\,dr}{s^{\frac {n-2}{n-1}} }
=
C \sup _{x \in \partial \o} \|\mathcal B \|_{L^{n-1, \infty}(\partial \Omega \cap B_r(x))}
\end{align}
for some constant $C=C(L_\Omega,  d_\o, n)$, for every  $x \in \partial \Omega$ and $r >0$.
An application of inequality \eqref{main6} with $v = \xi \,a(|\nabla u|) u_{x_i}$ , for $i=1, \dots n$, yields, via \eqref{main9},
\begin{multline}\label{main7}
 \int _{\partial \Omega } \xi^2 \, a(|\nabla u|)^2 |\nabla u|^2 |\mathcal B| \,d\mathcal H^{n-1}(x)  \\ \leq C \sup _{x \in \partial \o}  \|\mathcal B \|_{L^{n-1, \infty}(\partial \Omega \cap B_r(x))} \bigg(\int _\Omega  \xi ^2 a(|\nabla u|)^2 |\nabla ^2 u|^2\, dx  + \int _\Omega  |\nabla \xi |^2 a(|\nabla u|)^2 |\nabla u|^2\, dx \bigg)
\end{multline}
for some constant $C=C(L_\o, d_\o, n, s_a)$. Note that here we have made use of  the second inequality in \eqref{infsup} to infer that
\begin{equation}\label{feb20'}
|\nabla (a(|\nabla u|)u_{x_i})| \leq C \, a(|\nabla u|)|\nabla ^2u| \qquad 	\hbox{in $\o$,}
\end{equation}
for $i=1, \dots , n$, and for some constant $C=C(n, s_a)$.
Combining equations  \eqref{main3bis} and \eqref{main7} tells us that 
\begin{multline}\label{main10}
\Big[C_1(1- \varepsilon) - C_2 \sup _{x \in \partial \o} \|\mathcal B \|_{L^{n-1, \infty}(\partial \Omega \cap B_r(x)} \Big] \int _\Omega \xi ^2 a(|\nabla u|)^2 |\nabla ^2 u|^2\, dx \\ \leq \int _\Omega \xi^2 f^2\, dx 
+ \Big[C_2 \sup _{x \in \partial \o} \|\mathcal B \|_{L^{n-1, \infty}(\partial \Omega \cap B_r(x)} + \frac {C_3}\varepsilon \Big] \int _\Omega  |\nabla \xi |^2 a(|\nabla u|)^2 |\nabla u|^2\, dx
\end{multline}
for some constants $C_1=C_1(n, i_a)$, $C_2=C_2(L_\o, d_\o, n, s_a)$ and $C_3=C_3(n)$. 
 If condition \eqref{smalln} is fulfilled with $c= \frac {C_1}{C_2}$, then
  there exists $r_0>0$ such that 
$$C_1(1- \varepsilon) - C_2 \sup _{x \in \partial \o} \|\mathcal B \|_{L^{n-1, \infty}(\partial \Omega \cap B_r(x)}>0$$
%
%
if $0<r\leq r_0$ and $\varepsilon$ is sufficiently small. Therefore, by inequality \eqref{main10},
\begin{align}\label{main11}
\int _\Omega \xi ^2 a(|\nabla u|)^2 |\nabla ^2 u|^2\, dx & \leq C \int _\Omega \xi^2 f^2\, dx 
+  C \int _\Omega  |\nabla \xi |^2 a(|\nabla u|)^2 |\nabla u|^2\, dx
\end{align}
for some constant $C=C(L_\o, d_\o, n, i_a, s_a)$, if  $0<r\leq r_0$ in \eqref{feb22}.
\\ \color{black} In the case when $n=2$,   define
\begin{equation}\label{main5bis}
Q_1\textcolor{black}{(r)} = \sup _{x \in \partial \o} \sup _{E \subset \partial \Omega\textcolor{black}{\cap B_r(x)}} \frac{\int _E |\mathcal B| \,d \mathcal H^1 (y)}{{\rm cap}_{B_1(x)}(E)}\, \quad \quad \hbox{for $r \in (0,1)$,}
\end{equation}
where ${\rm cap}_{B_1(x)}(E)$ stands for the  capacity  of the set
 $E$  given by
\begin{equation}\label{cap2}
{\rm cap}_{B_1(x)}(E) = \inf \bigg\{\int _{B_1(x)} |\nabla v|^2\, dy : v\in C^1_0(B_1(x)), v\geq 1 \,\hbox{on }\, E\bigg\}.
\end{equation}
A counterpart of inequality \eqref{main6}  reads
\begin{equation}\label{main6bis}
\int _{\partial \Omega \cap B_{\textcolor{black}{r}}(x)} v^2 \, |\mathcal B| \, d \mathcal H^1 (y) \leq C \textcolor{black}{Q_1(r)}
\int _{\Omega \cap B_{\textcolor{black}{r}}(x)}|\nabla v|^2 \, dy
\end{equation}
for every $x \in \partial \Omega$, $r \in (0, 1)$ and $v \in C^1_0(B_{\textcolor{black}{r}}(x))$, where $C=C(L_\o, d_\o)$. 
\\ A borderline version of the trace inequality -- see e.g.  \cite[Section 7.6.4]{AH} -- \color{black} ensures that there exists a constant $C=C(L_\o, d_\o, n)$ such that 
\begin{equation}\label{trace2}
\sup _{E \subset \partial \Omega \cap B_1(x)}\frac{\big(\frac 1{\mathcal H^1 (E)}\int _E v \,d \mathcal H^1 (y)\big)^2}{\log \big(1+ \frac{\mathcal H^1(\partial \Omega \cap B_1(x))}{\mathcal H^1(E)}\big)} \leq C \int _{\Omega \cap B_1(x)}|\nabla v|^2 \, dy
%
\end{equation}
for every  $x \in \partial \Omega$, and \textcolor{black}{$v \in C^1_0(B_1(x))$}. Notice that the left-hand side of \eqref{trace2} is equivalent to the norm in an Orlicz space associated with the Young function $e^{t^2}-1$.
The choice of trial functions $v$ in \eqref{trace2} such that $v\geq 1$ on $E$ yields, 
via definition \eqref{cap2},
\begin{equation}\label{main8bis}
\frac 1{\log \big(1 + \frac{C}{\mathcal H^1 (E)}\big)} \leq C {\rm cap}_{B_1(x)} (E)\,,
\end{equation}
for some constant $C=C(L_\o , d_\o)$, and for every set $E \subset \partial \Omega \cap B_1(x)$. 
Thanks to \eqref{main8bis} and to the Hardy-Littlewood inequality again, 
\begin{align}\label{feb20}
Q_1(r) &\leq \sup _{x \in\partial \Omega} \sup _{E \subset \partial \Omega\textcolor{black}{\cap B_r(x)}} \frac{\int _0^{\mathcal H^1(E)} (|\mathcal B| _{|\partial \Omega\textcolor{black}{\cap B_r(x)}})^*(r)\, dr}{{\rm cap}_{B_1(x)}(E)}
\\ \nonumber & \leq C \sup _{x \in\partial \Omega} \sup_{s \in (0, \mathcal H^1( \partial \Omega\textcolor{black}{\cap B_r(x)}))} \log \Big(1 + \frac{C}{s}\Big)\int _0^{s} (|\mathcal B| _{|\partial \Omega\textcolor{black}{\cap B_r(x)}})^*(r)\, dr 
\\ \nonumber & = C  \sup _{x \in\partial \Omega}\|\mathcal B \|_{L^{1, \infty}\log L (\partial \Omega\textcolor{black}{\cap B_r(x)})}
\end{align}
for some constant $C=C(L_\o , d_\o)$, and for $r \in (0,1)$.
%
%
%
On exploiting \eqref{feb20} instead of \eqref{main9},   and arguing as  in the case when $n \geq 3$, yield \eqref{main11} also for $n=2$.
\\ \color{black} When $B_r(x) \subset \subset \o$, the derivation of \eqref{main11} is even simpler, and follows directly from \eqref{main1}, \eqref{main2} and \eqref{main4}, since the boundary integral on the rightmost side    of \eqref{main2} vanishes in this case.
\\ Now, let  $\{B_{r_k}\}_{k\in K}$ be a finite covering of $\overline \o$ by balls $B_{r_k}$, with  $r_k \leq r_0$, such that either  $B_{r_k}$ is centered on $\partial \o$, or $B_{r_k}\subset  \subset \o$. Note that this covering can be chosen in such a way that the \textcolor{black}{multiplicity of overlapping of the balls $B_{r_k}$ only depends on $n$.} Let 
$\{\xi _k\}_{k\in K}$ be a family of functions such that $\xi_k \in C^\infty _0(B_{r_k})$ and 
$\{\xi_k^2\}_{k\in K}$ is a partition of unity associated with the covering $\{B_{r_k}\}_{k\in K}$.
Thus  $\sum _{k\in K} \xi _k^2 = 1$ in $\overline \o$.  On applying inequality \eqref{main11} with $\xi = \xi _k$ for each $k$, and adding the resulting inequalities one obtains that
\begin{align}\label{main12}
\int _\Omega  a(|\nabla u|)^2 |\nabla ^2 u|^2\, dx & \leq C \int _\Omega  f^2\, dx 
+  C \int _\Omega  a(|\nabla u|)^2 |\nabla u|^2\, dx
\end{align}
for some constant $C= C(L_\o, d_\o, n, i_a, s_a)$. 
\\ A version of the Sobolev inequalty  entails that, for every $\sigma >0$, there exists a constant 
$C=C(L_\o, d_\o,   n, \sigma)$ 
such that
\begin{align}\label{main13}
\int _\Omega v^2\, dx \leq \sigma \int _\Omega |\nabla v|^2\, dx + C \bigg(\int _\Omega |v|\, dx \bigg)^2
\end{align}
for every $v \in W^{1,2}(\Omega)$ (see e.g. \cite[Proof of Theorem 1.4.6/1]{Mabook}). Applying inequality \eqref{main13} with $v = a(|\nabla u|) u_{x_i}$, $i=1, \dots , n$, an recalling \eqref{feb20'} tell  us that
\begin{align}\label{main14}
\int _\Omega  a(|\nabla u|)^2 |\nabla u|^2\, dx \leq \sigma C_1 \int _\Omega a(|\nabla u|)^2 |\nabla ^2 u|^2\, dx + C_2 \bigg(\int _\Omega a(|\nabla u|)|\nabla u|\, dx \bigg)^2
\end{align}
for some constant $C_1=C_1(n, s_a)$ and $C_2= C_2(L_\o, d_\o, n,  s_a, \sigma )$.
%
%
%
On choosing $\sigma = \tfrac{1}{2CC_1}$, where $C$ is the constant appearing in \eqref{main12}, and  combining inequalities \eqref{main12}, \eqref{main14} and \eqref{estgraddir} we conclude that
\begin{align}\label{main16}
\int _\Omega  a(|\nabla u|)^2 |\nabla ^2u|^2\, dx \leq C \int _\Omega  f^2\, dx 
\end{align}
for some constant $C=C(L_\o, d_\o, n, i_a , s_a)$. 
Inequalities \eqref{main14}, \eqref{main16} and \eqref{estgraddir} imply, via \eqref{feb20'}, that
\begin{equation}\label{main100aus}
\|a(|\nabla u|) \nabla u\|_{W^{1,2}(\o)} \leq C \|f\|_{L^2(\o)}
\end{equation}
for some constant $C=C(L_\o, d_\o, n, i_a, s_a)$. In particular, the dependence of the constant $C$ in \eqref{main100aus}  is in fact just through an upper bound for the quantities $L_\o, d_\o,  s_a$, and through a lower bound for $i_a$. This is crucial in view of the next steps.

\smallskip
\par\noindent
 \emph{Step 2}. Here we remove assumptions \eqref{abound} and \eqref{Asmooth}.  To this purpose, we make use of a family of functions $\{a_\ep\}_{\ep \in (0,1)}$, with 
 $a_\varepsilon : [0, \infty ) \to
(0, \infty )$, satisfying the following properties:
\begin{equation}\label{aeps}
a_\varepsilon : [0, \infty) \to [0, \infty) \quad \hbox{and}  \quad \varepsilon \leq    a_\varepsilon (t)   \leq \varepsilon ^{-1} \quad
\hbox{for $t \geq 0$;}
\end{equation}
\begin{equation}\label{indici}
\min \{i_a , 0\} \leq i_{a_\varepsilon} \leq s_{a_\varepsilon} \leq
\max \{s_a , 0\};
\end{equation}
\begin{equation}\label{conva}
 \lim _{\varepsilon \to 0} a_\varepsilon (|\xi |) \xi = a (|\xi |) \xi \qquad
\hbox{uniformly in $\{\xi \in \R ^{n} : |\xi | \leq M\}$  for every
$M>0$;}
\end{equation}
\color{black}
the function $\mathcal A_\varepsilon : \rn \to [0, \infty)$, defined as $\mathcal A_\varepsilon(\eta) = a_\varepsilon(|\eta|)$ for $\eta \in \rn$, is such that 
\begin{equation}\label{Aepssmooth}
\mathcal A_\varepsilon  \in C^\infty(\rn).
\end{equation}
\color{black}
The construction of a family of functions enjoying these properties can be accomplished on combining 
 \cite[Lemma 3.3]{cmCPDE} and  \cite[Lemma 4.5]{cmARMA}.  
Now, let $u_\varepsilon$ be the solution to the problem 
\begin{equation}\label{eqep}
\begin{cases}
- {\rm div} (a_\ep(|\nabla u_\ep|)\nabla u_\ep ) = f(x)  & {\rm in}\,\,\, \o \\
 u_\ep =0  &
{\rm on}\,\,\,
\partial \o \,.
\end{cases}
\end{equation}
Owing to  \eqref{aeps} and \eqref{Aepssmooth}, the assumptions of Step 1 are fulfilled by problem \eqref{eqep}. Thus, as a consequence of \eqref{main100aus}, there exists    a constant $C= C(L_\o, d_\o, n, i_a, s_a)$ such that 
\begin{equation}\label{main17}
\|a_\varepsilon(|\nabla u_\varepsilon|) \nabla u_\varepsilon\|_{W^{1,2}(\o)} \leq C \|f\|_{L^2(\o)}
\end{equation}
%
%
for $\varepsilon \in (0, 1)$.
Observe that the constant  $C$ in \eqref{main17} is actually  independent of $\ep$, thanks to \eqref{indici}. 
By \eqref{main17}, there exists a sequence $\{\ep _k\}$ and a function $U : \Omega \to \rn$  such that $U \in W^{1,2}(\Omega)$, 
\begin{equation}\label{main19}
a_{\ep _k}(|\nabla u_{\ep _k}|) \nabla u_{\ep _k} \to U\quad \hbox{in $L^2(\Omega)$} \quad \hbox{and} \quad a_{\ep _k}(|\nabla u_{\ep _k}|) \nabla u_{\ep _k} \rightharpoonup U\quad \hbox{in $W^{1,2}(\Omega)$,}
\end{equation}
where the arrow $\lq\lq \rightharpoonup  "$ stands for weak convergence.
On the other hand,   
a global estimate for $\|u_{\ep _k} \|_{L^\infty (\o)}$ following from a result of \cite{Talenti_orlicz}, coupled with 
a local gradient estimate of \cite[Theorem 1.7]{Li} ensures that $u_{\ep _k} \in C^{1,\alpha}_{\rm loc}(\Omega)$, and  that for any open set $\o ' \subset \subset \Omega$ there exists a constant $C$ such that  
\begin{align}\label{main21}
\|u_{\ep _k}\|_{C^{1,\alpha}(\Omega ')} \leq C\,
\end{align}
for $k \in \mathbb N$.
%
%
%
%
Thus,   there exists a function $v \in C^{1}(\Omega)$ such that, on taking, if necessary,   a subsequence, 
\begin{equation}\label{main22}
u_{\ep _k} \to v\, \quad \hbox{and}\quad \nabla u_{\ep _k} \to \nabla v \quad \hbox{pointwise in $\Omega$.}
\end{equation}
In particular, 
\begin{equation}\label{main25'}
a(|\nabla v|) \nabla v = U,
\end{equation}
and hence
\begin{equation}\label{main25}
a(|\nabla u|) \nabla u  \in W^{1,2}(\Omega)\,.
\end{equation}
Testing the equation in \eqref{eqep} with any function $\varphi \in C^\infty_0(\Omega)$ yields
\begin{equation}\label{main23}
\int _\Omega a_{\ep_k}(|\nabla u_{\ep_k}|) \nabla u_{\ep_k} \cdot \nabla \varphi \, dx = \int_\Omega f \, \varphi \, dx\,.
\end{equation}
Owing to \eqref{main19} and \eqref{main25'},  on passing to the limit in \eqref{main23} as $k \to \infty$ one deduces that
%
\begin{equation}\label{main24}
\int _\Omega a(|\nabla v|) \nabla v \cdot \nabla \varphi \, dx = \int_\Omega f \, \varphi \, dx\,.
\end{equation}
Thus $v=u$, the weak solution to problem \eqref{eqdirichlet}. Furthermore, by \eqref{main17}, we obtain via \eqref{main19} and \eqref{main25'}  that
\begin{equation}\label{main27}
\|a(|\nabla u|) \nabla u\|_{W^{1,2}(\o)} \leq C \|f\|_{L^2(\o)}
\end{equation}
%
%
for some constant $C= C(L_\o, d_\o,  n, i_a, s_a)$.

\smallskip
\par\noindent
\emph{Step 3}. Here, we remove assumption \eqref{omegasmooth}. 
Via  smooth approximation of the functions which locally describe $\partial \Omega$, one can construct a sequence
 $\{\Omega _m\}$  of open sets in $\rn$ such that $\partial \Omega _m \in C^\infty$, $\Omega \subset \Omega _m$, $\lim _{m \to \infty}|\Omega _m \setminus \Omega| = 0$, and the Hausdorff distance between $\Omega _m$ and $\Omega$ tends to $0$ as $m \to \infty$. Also, there exists a constant $C=C(\Omega)$ such that
\begin{equation}\label{feb35}
L_{\o _m} \leq C L_\o \quad \hbox{and} \quad d_{\o _m} \leq C d_\o
\end{equation}
for $m \in \mathbb N$. Moreover, 
 although  smooth functions are   neither dense in $W^2L^{n-1, \infty}$  if $n\geq 3$, nor in $W^2L^{1, \infty}\log L$ if $n=2$, one has that
$$
\sup _{x \in \partial \o}\|\mathcal B _m\|_{L^{n-1, \infty}(\partial \Omega _m \cap B_r(x))}
\leq C \sup _{x \in \partial \o}\|\mathcal B \|_{L^{n-1, \infty}(\partial \Omega \cap B_r(x))} \quad \hbox{if $n \geq 3$,}$$
or
$$
\sup _{x \in \partial \o}\|\mathcal B _m\|_{L^{1, \infty}\log L(\partial \Omega _m \cap B_r(x))}
\leq C \sup _{x \in \partial \o}\|\mathcal B \|_{L^{1, \infty}\log L(\partial \Omega \cap B_r(x))} \quad \hbox{if $n =2$,}$$
for some
constant $C=C(\Omega)$, where $\mathcal B _m$ denotes the second fundamental form on  $\partial \Omega _m$. 
\\ Let $u_m$ be the weak solution to the Dirichlet problem 
\begin{equation}\label{eqm}
\begin{cases}
- {\rm div} (a(|\nabla u_m|)\nabla u_m ) = f(x)  & {\rm in}\,\,\, \o _m \\
 u_m =0  &
{\rm on}\,\,\,
\partial \o _m \,,
\end{cases}
\end{equation}
where $f$ still fulfils \eqref{fsmooth}, and  is extended by $0$ outside $\Omega$. By inequality \eqref{main27} of Step 2,  
\begin{equation}\label{main27aus}
\|a(|\nabla u_m|) \nabla u_m\|_{W^{1,2}(\o_m)} \leq C \|f\|_{L^2(\o_m)}= C \|f\|_{L^2(\o)},
\end{equation}
%
%
the constant $C$ being independent of $m$, by the properties of $\o_m$ mentioned above. 
\\
Thanks to \eqref{main27aus}, the sequence  
$\{a(|\nabla u_m |) \nabla u_m\}$ is  bounded in $W^{1,2}(\Omega)$,
and hence  there exists a subsequence, still denoted by $\{u_m\}$ and a function $U : \Omega \to \rn$ such that $U \in W^{1,2}(\Omega)$,
\begin{equation}\label{main32}
a(|\nabla u_{m}|) \nabla u_{m} \to U\quad \hbox{in $L^2(\Omega )$} \quad \hbox{and} \quad a(|\nabla u_{m}|) \nabla u_{m} \rightharpoonup U\quad \hbox{in $W^{1,2}(\Omega )$}.
\end{equation}
By the local gradient estimate recalled in Step 2, there exists $\alpha \in (0,1)$ such that $u _m\in C^{1,\alpha}_{\rm loc}(\Omega)$, and for every open set $\o ' \subset \subset \Omega$ there exists a constant $C$, independent of $m$, such that  
\begin{align}\label{main29}
\|u_{m}\|_{C^{1,\alpha}(\Omega ')} \leq C\,.
\end{align}
Thus,  on taking, if necessary,  a further subsequence,  
\begin{equation}\label{main30}
u_{m} \to v\, \quad \hbox{and}\quad \nabla u_{m} \to \nabla v \quad \hbox{pointwise in $\Omega$,}
\end{equation}
for some function $v \in C^1(\o)$. 
In particular,
\begin{equation}\label{main31}
a(|\nabla u_{m}|) \nabla u_{m} \to a(|\nabla v|) \nabla v \quad \hbox{pointwise in $\Omega$.}
\end{equation}
By \eqref{main31} and \eqref{main32}, 
\begin{equation}\label{main100}
a(|\nabla v|) \nabla v   = U \in W^{1,2}(\Omega )\,.
\end{equation}
 Given any function $\varphi \in C^\infty _0(\o)$, on passing to the limit as $m \to \infty$ in the weak formulation of problem \eqref{eqm}, namely in the equation
\begin{equation}\label{main101}
\int _{\Omega_m} a(|\nabla u_m|) \nabla u_m \cdot \nabla \varphi \, dx = \int_{\Omega _m} f \, \varphi \, dx\,,
\end{equation}
we infer from \eqref{main32} and \eqref{main100} that
$$
\int _\Omega a(|\nabla v|) \nabla v \cdot \nabla \varphi \, dx = \int_\Omega f \, \varphi \, dx\,.
$$
Therefore, $u=v$, the weak solution to problem \eqref{eqdirichlet}. 
Furthermore, owing to \eqref{main27aus}, \eqref{main32} and \eqref{feb20'},
\begin{equation}\label{main35}
\|a(|\nabla u|) \nabla u\|_{W^{1,2}(\o)} \leq C \|f\|_{L^2(\o)}
\end{equation}
%
%
for some constant $C=C(L_\o, d_\o, n, i_a, s_a)$.

\smallskip
\par\noindent
\emph{Step 4}. We conclude by removing the remaining additional assumption \eqref{fsmooth}. Let $f \in L^2(\Omega)$. 
 Owing to \eqref{gen1}, given  any sequence $\{f_k\} \subset C^\infty_0(\Omega )$ such that    $f_k \to f$ in $L^2(\o)$,  the sequence    $\{u_k\}$  of the weak solutions to the Dirichlet problems 
\begin{equation}\label{eqk}
\begin{cases}
- {\rm div} (a(|\nabla u_k|)\nabla u_k ) = f_k  & {\rm in}\,\,\, \o  \\
 u_k =0  &
{\rm on}\,\,\,
\partial \o  \,,
\end{cases}
\end{equation} 
fullfils
\begin{equation}\label{main36}
u_k \to u \quad \hbox{and} \quad \nabla u_k \to \nabla u \quad \hbox{a.e. in $\o$.}
\end{equation}
%
By inequality \eqref{main35} of the previous step, we have that $a(|\nabla u_k|) \nabla u_k \in W^{1,2}(\Omega)$, and there exist  constants $C_1$ and $C_2$, independent of $k$, such that
\begin{align}\label{main40}
\|a(|\nabla u_k|) \nabla u_k\|_{W^{1,2}(\o)} \leq C_1 \|f_k\|_{L^2(\o)} \leq C_2 \|f\|_{L^2(\o)}\,.
\end{align}
Hence, the sequence $\{a(|\nabla u_k|)\nabla u_k\}$ is uniformly bounded in $ W^{1,2}(\Omega)$, and there exists a subsequence, still indexed by $k$, and a function $U :  \Omega \to \mathbb R^n$ such that $U \in W^{1,2}(\Omega)$ and 
\begin{equation}\label{main38}
a_{k}(|\nabla u_{k}|) \nabla u_{k} \to U\quad \hbox{in $L^2(\Omega)$} \quad \hbox{and} \quad  a_{k}(|\nabla u_{k}|) \nabla u_{k} \rightharpoonup U\quad \hbox{in $W^{1,2}(\Omega)$}.
\end{equation}
From \eqref{main36} we thus infer that $a(|\nabla u|) \nabla u =U  \in W^{1,2}(\Omega)$,
and the second inequality in \eqref{seconddir2} follows via \eqref{main40} and     \eqref{main38}.  The first inequality is easily verified, via \eqref{feb20'}.
%
%
The statement concerning the solution to  the Dirichlet problem \eqref{eqdirichlet} is thus fully proved. 

\smallskip
\par We point out hereafter the  changes required for  the solution to the Neumann problem \eqref{eqneumann}.
\\ \emph{Step 1}. The additional assumption \eqref{0mean} has to be coupled with \eqref{fsmooth}. Moreover, since $\frac {\partial u}{\partial \nu}=0$ on $\partial \o$, the middle term in the chain \eqref{main3} is replaced with 
$$
\bigg|- 
\int _{\partial \Omega } \xi^2 a(|\nabla u|)^2  \mathcal B (\nabla _T \, u , \nabla _T \, u)
 \,d\mathcal H^{n-1}(x)\bigg| \,.$$
\\ \emph{Step 2}. The Dirichlet boundary condition in problem \eqref{eqep} must, of course, be replaced with the Neumann condition $\frac{\partial u_\ep}{\partial \nu} =0$. The solution of the resulting Neumann problem is only unique up to additive constants. A bound of the form $\|u_{\ep_k}- c_k\|_{L^\infty (\o)}\leq C$ now holds for a suitable sequence  $\{c_k\}$ with $c_k\in \mathbb R$ \cite{CianchiPRSE}. Hence, $u_{\ep _k}$ has to be replaced with $u_{\ep_k}- c_k$ in equations  \eqref{main21} and \eqref{main22}. Moreover, the test functions $\varphi$ in equation \eqref{main23} now belong to $W^{1, \infty}(\o)$.

\smallskip
\par\noindent
 \emph{Step 3}. 
The Dirichlet problem \eqref{eqm} has to be replaced with the Neumann problem with boundary condition $\frac{\partial u_m}{\partial \nu} =0$. Accordingly, the corresponding sequence of solutions $\{u_m\}$ has  to be normalized by a suitable sequence of additive constants.
\\ Passage to the limit as $m \to \infty$ in equation \eqref{main101}  can be justified as follows. Extend any  test function $\varphi \in W^{1, \infty}(\o)$  to a function in $W^{1, \infty}(\mathbb R^n)$, still denoted by $\varphi$. The left-hand side of equation \eqref{main101} can be split as 
\begin{equation}\label{main101'}
\int _{\o _m} a(|\nabla u_m|)\nabla u_m \cdot \nabla \varphi \, dx = 
\int _{\o } a(|\nabla u_m|)\nabla u_m \cdot \nabla \varphi \, dx + \int _{\o _m \setminus \o} a(|\nabla u_m|)\nabla u_m \cdot \nabla \varphi \, dx\,.
\end{equation}
The first integral on the right-hand side of \eqref{main101'} converges to 
$$\int _{\o} a(|\nabla v|)\nabla v \cdot \nabla \varphi \, dx\,$$
as $m \to \infty$, owing to \eqref{main32} and \eqref{main100}. The second integral tends to $0$, by \eqref{main27aus} and the fact that $|\o _m \setminus \o| \to 0$. 

\smallskip
\par\noindent
\emph{Step 4}. The sequence of approximating functions $\{f_k\}$ has to fulfill the additional compatibility condition $\smallint _\o f_k(x)\, dx =0$ for $k \in \mathbb N$.  Moreover, the Dirichlet boundary condition in problem \eqref{eqk} has to be replaced with the Neumann condition $\frac {\partial u_k}{\partial \nu}=0$ on $\partial \Omega$.
\qed

\medskip
\par\noindent
{\bf Proof of Theorem \ref{secondconvex}}.  The proof parallels (and is even simpler than) that of Theorem \ref{seconddir}. We limit oureselves to pointing out the variants and simplifications needed.
\\
 \emph{Step 1}. Assume
that $\Omega$, $a$ and $f$ are as in Step 1 of the proof of
Theorem \ref{seconddir} and that, in addition,  $\Omega $ is convex.
One can proceed as in that proof, and
exploit the fact that  the right-hand side of equation \eqref{grisv1} is nonnegative owing to the convexity of $\o$, since it reduces to either
$$     - {\rm tr}\mathcal B \bigg(\frac {\partial u}{\partial \nu}\bigg)^2  \geq 0 \quad  \hbox{or}  \quad 
-\mathcal B (\nabla _T \,u , \nabla _T \,u) \geq 0 \qquad \quad
\hbox{on $\partial \Omega$ \,,}$$
according to whether $u$ is the solution to the Dirichlet problem \eqref{eqdirichlet}, or to the Neumann problem \eqref{eqneumann}. Therefore, inequality \eqref{main3bis} can be replaced with the stronger inequality
\begin{align}\label{main3ter}
C(1-\varepsilon) \int _\Omega \xi^2 a(|\nabla u|)^2 |\nabla^2 u|^2\,dx  \leq 
\int _\Omega \xi^2 f^2\, dx  + \frac {C'} \varepsilon \int _\Omega |\nabla \xi|^2 a(|\nabla u|)^2 |\nabla u|^2\,dx \,.
\end{align}
 Starting from this inequality, instead of \eqref{main3bis},
estimate \eqref{main100} follows analogously.

\smallskip
\par\noindent \emph{Step
2}. The proof is the same as that of Theorem \ref{seconddir}.

\smallskip
\par\noindent \emph{Step
3}. The proof is 
analogous to that of Theorem \ref{seconddir}, save
that the approximating domains $\Omega _m$ have to be chosen in such a way that they are
convex.

\smallskip
\par\noindent \emph{Step
4}. The proof is the same as that of Theorem \ref{seconddir}.
\qed

\section{Local estimates}\label{proofsloc}

Here, we provide a proof of Theorem \ref{secondloc}. The generalized local solutions to equation \eqref{localeq} considered in the statement can be defined as follows.
\par Assume that  $f \in L^q_{\rm loc}(\Omega)$ for some $q\geq1$.  
A function $u \in \mathcal T ^{1,1}_{\rm loc}(\o)$ is called a  generalized local solution to equation \eqref{localeq} if $a(|\nabla u|)\nabla u\in L^1_{\rm loc}(\o)$, equation \eqref{231} holds for every   $\varphi \in C^\infty _0(\o)$, and
there exists a sequence $\{f_k\}\subset  C^\infty_0 (\o)$ and a correpsonding sequence of local weak solutions $\{u_k\}$ to equation \eqref{localeq}, with $f$ replaced by $f_k$,
such that 
$f_k \to f$   in $L^q (\o')$,
\begin{equation}\label{approxuloc}
u_k \to u \quad \hbox{and} \quad \nabla u_k \to \nabla u \quad \hbox{a.e. in $\o$,}
\end{equation}
and
\begin{equation}\label{approxaloc}
\lim _{k \to \infty} \int_{\o'} a(|\nabla u_k|) |\nabla u_k| \, dx =\int_{\o'} a(|\nabla u|) |\nabla  u|\, dx\,
\end{equation}
for every open set $\o' \subset \subset \o$.
\par Note that, by the results from \cite{cmapprox} recalled at the beginning of Section \ref{proofs}, the generalized solutions to the boundary value problems \eqref{eqdirichlet} and \eqref{eqneumann} are, in particular, generalized local solutions to equation \eqref{localeq}.

%
%

\smallskip
\par\noindent
{\bf Proof of Theorem \ref{secondloc}}.  This proof  follows the outline  of that of Theorem \ref{seconddir}. Some variants are however required, due to the local nature of the result. Of course, the step concerning the approximation of $\Omega$ by domains with a smooth boundary  is not needed at all.

\smallskip
\par\noindent
\emph{Step 1}. Assume the additional conditions \eqref{fsmooth} on $f$, and \eqref{abound} -- \eqref{Asmooth} on $a$, and let $u$ be a local weak solution to equation \eqref{localeq}. Thanks to the current assumption on $a$ and $f$, the function $u$ is in fact a classical smooth solution. Let $B_{2R}$ be any ball such that $B_{2R} \subset \subset \Omega$, and let $R\leq \sigma < \tau \leq 2R$. An application of inequality \eqref{main3bis}, with $\varepsilon = \tfrac 12$ and any function $\xi \in C^{\infty}_0(B_\tau)$ such that $\xi =1 $ in $B_\sigma $ and $|\nabla \xi|\leq C/ (\tau -\sigma)$ for some constant $C=C(n)$, tells us that
\begin{align}\label{loc10}
\int _{B_\sigma} a(|\nabla u|)^2 |\nabla^2 u|^2\,dx   \leq 
C \int _{B_{2R}} f^2\, dx  + \frac C{(\tau - \sigma)^2}\int _{B_\tau \setminus B_\sigma} a(|\nabla u|)^2 |\nabla u|^2\,dx 
\end{align}
for some constant $C=C(n, i_a, s_a)$. We claim that there exists a constant $C=C(n)$ such that
\begin{equation}\label{sobolevring}
\int _{B_\tau \setminus B_\sigma}  v^2\, dx \leq 	\frac {\delta ^2}{(\tau - \sigma)^2} \int _{B_\tau \setminus B_\sigma}  |\nabla v|^2\, dx + \frac{C(\tau - \sigma) R^{n-1}}{ \delta^{n}} \bigg(\int _{B_\tau \setminus B_\sigma}   |v|\, dx\bigg)^2
\end{equation}
for every $\delta >0$ and  every $v \in W^{1,2}(B_\tau \setminus B_\sigma)$, provided  that $R$, $\tau$ and $\sigma$ are as above. This claim can be verified as follows. Denote by $Q_r$ a cube of sidelength $r>0$.
The inequality 
\begin{equation}\label{sobring1}
\int _{Q_1}v^2\, dx \leq C_1 \int _{Q_1}|\nabla v|^2\, dx + C_2 \bigg(\int _{Q_1} |v|\, dx \bigg)^2
\end{equation}
holds for every $v \in W^{1,2}(Q_1)$, for suitable constants $C_1=C_1(n)$ and $C_2(n)$. Given $\ep >0$, a scaling argument tells us that a parallel inequality holds in  $Q_\ep$, with $C_1$ replaced with $C_1 \ep ^2$ and $C_2$ replaced with $C_2\ep^{-n}$. A covering argument for $Q_1$ by cubes of sidelength $\ep$ then yields inequality \eqref{sobring1} with $C_1$ and $C_2$ replaced by $C_1 \ep ^2$ and $C_2\ep^{-n}$, respectively. Another scaling argument, applied to the resulting inequality in $Q_1$, provides us with the inequality
\begin{equation}\label{sobring2}
\int _{Q_\delta}v^2\, dx \leq C_1 (\ep \delta)^2\int _{Q_\delta}|\nabla v|^2\, dx + C_2(\ep \delta)^{-n} \bigg(\int _{Q_\delta} |v|\, dx \bigg)^2
\end{equation}
for every $v \in W^{1,2}(Q_\delta)$. Via a covering argument for $B_2\setminus B_1$ by (quasi)-cubes of suitable sidelength $\delta$, one  infers from \eqref{sobring2} that 
\begin{equation}\label{sobring3}
\int _{B_2 \setminus B_1}v^2\, dx \leq C \ep^2\int _{B_2 \setminus B_1}|\nabla v|^2\, dx + C \ep ^{-n} \bigg(\int _{B_2 \setminus B_1} |v|\, dx \bigg)^2
\end{equation}
for a suitable constant $C=C(n)$. Inequality \eqref{sobolevring} can be derived from \eqref{sobring3} on mapping $B_2 \setminus B_1$ into  $B_\tau \setminus B_\sigma$ via   the bijective map $\Phi : B_2 \setminus B_1 \to B_\tau \setminus B_\sigma$ defined as  
$$\Phi (x) = \frac {x}{|x|} \big[\sigma + (|x|-1)(\tau - \sigma )\big] \quad \hbox{for $x \in B_2 \setminus B_1 $,}$$
and making use of the fact that 
$$c_1 (\tau - \sigma) R^{n-1} \leq |\det (\nabla \Phi (x))|\leq  c_2(\tau - \sigma) R^{n-1} \quad \hbox{for $x \in B_2 \setminus B_1 $}$$
and 
$$|\nabla (\Phi^{-1}) (y) | \geq c_1 (\tau - \sigma) \quad \hbox{for $y \in B_\tau \setminus B_\sigma $,}$$
for suitable positive constants $c_1=c_1(n)$ and $c_2=c_2(n)$.
\\
Choosing $\delta = (\tau -\sigma)^2$ in inequality \eqref{sobolevring}, and applying the resulting inequality with $v= a(|\nabla u|)u_{x_i}$, for $i=1 \dots ,n$ yields
\begin{multline}\label{sobolevring1}
\frac {1}{(\tau - \sigma)^2}\int _{B_\tau \setminus B_\sigma}  a(|\nabla u|)^2 |\nabla u|^2\, dx \\ \leq 	C  \int _{B_\tau \setminus B_\sigma}  a(|\nabla u|)^2 |\nabla^2 u|^2\, dx + \frac{CR^{n-1}}{ (\tau - \sigma) ^{2n-1}} \bigg(\int _{B_\tau \setminus B_\sigma}   a(|\nabla u|)|\nabla u|\, dx\bigg)^2
\end{multline}
for some constant $C=C(n, s_a)$. Observe that in \eqref{sobolevring1} we have also made use of equation \eqref{feb20'}. Inequalities \eqref{loc10} and \eqref{sobolevring1} imply that
\begin{align}\label{loc11}
\int _{B_\sigma} a(|\nabla u|)^2 |\nabla^2 u|^2\,dx  &  \leq 
 C  \int _{B_\tau \setminus B_\sigma}  a(|\nabla u|)^2 |\nabla^2 u|^2\, dx 
\\ \nonumber & \quad 
+
C \int _{B_{2R}} f^2\, dx  
+ \frac{CR^{n-1}}{ (\tau - \sigma) ^{2n-1}} \bigg(\int _{B_{2R}}   a(|\nabla u|)|\nabla u|\, dx\bigg)^2
\end{align}
for some constant $C=C(n, i_a, s_a)$. Adding the quantity $C\int _{B_\sigma} a(|\nabla u|)^2 |\nabla^2 u|^2\,dx$ to both sides of inequality \eqref{loc11}, and dividing through the resulting inequality by $(1+C)$ enable us to deduce that
\begin{align}\label{loc12}
\int _{B_\sigma} a(|\nabla u|)^2 |\nabla^2 u|^2\,dx  &  \leq 
 \frac {C}{1+C}  \int _{B_\tau }  a(|\nabla u|)^2 |\nabla^2 u|^2\, dx 
\\ \nonumber & \quad 
+
C' \int _{B_{2R}} f^2\, dx  
+ \frac{C'R^{n-1}}{ (\tau - \sigma) ^{2n-1}} \bigg(\int _{B_{2R}} a(|\nabla u|)|\nabla u|\, dx\bigg)^2
\end{align}
for positive constants $C=C(n, i_a, s_a)$ and $C'=C'(n, i_a, s_a)$. Inequality \eqref{loc12}, via a standard iteration argument (see e.g. \cite[Lemma 3.1, Chapter 5]{giaquinta}), entails that
\begin{align}\label{loc13}
\int _{B_R} a(|\nabla u|)^2 |\nabla^2 u|^2\,dx   
\leq C \int _{B_{2R}} f^2\, dx  
+ \frac{C}{R^{n}} \bigg(\int _{B_{2R}} a(|\nabla u|)|\nabla u|\, dx\bigg)^2
\end{align}
for some constant $C=C(n, i_a, s_a)$. On the other hand, a scaling argument applied to the Sobolev inequality \eqref{main13}, with $\Omega = B_1$ and $\sigma =1$, tells us that there exists a constant $C=C(n, s_a)$ such that 
 \begin{align}\label{loc14}
\int _{B_R} a(|\nabla u|)^2 |\nabla u|^2\,dx  &  \leq 
  \int _{B_R}  a(|\nabla u|)^2 |\nabla^2 u|^2\, dx 
+ \frac{C}{ R^{n}} \bigg(\int _{B_{R}}   a(|\nabla u|)|\nabla u|\, dx\bigg)^2\,.
\end{align}
Coupling inequality \eqref{loc13} with \eqref{loc14} yields
\begin{align}\label{loc15}
\|a(|\nabla u|)\nabla u\|_{W^{1,2}(B_R)} \leq C\big(\|f\|_{L^2(B_{2R})} + R^{-\frac n2}\|a(|\nabla u|)\nabla u\|_{L^1(B_{2R})}\big)
\end{align}
for some constant $C=C(n, i_a, s_a)$.

\smallskip
\par\noindent
\emph{Step 2}. Assume that $u$ is a local solution to equation \eqref{localeq}, with $a$ as in the statement, and $f$ still fulfilling \eqref{fsmooth}. 
One has that $u \in L^\infty_{\rm loc} (\Omega)$. This follows from \cite[Theorem 5.1]{Kor}, or 
from  gradient regularity results of \cite{Baroni} or \cite{DKS}. As a consequence, 
by \cite[Theorem 1.7]{Li}, $u \in C^{1,\alpha}_{\rm loc}(\Omega)$ for some $\alpha \in (0,1)$. 
Next, consider   
a family of functions $\{a_\ep\}_{\ep \in (0,1)}$ satisfying properties  \eqref{aeps} -- \eqref{Aepssmooth}. 
Denote by $u_\varepsilon$ the solution to the problem 
\begin{equation}\label{eqeploc}
\begin{cases}
- {\rm div} (a_\ep(|\nabla u_\ep|)\nabla u_\ep ) = f(x)  & {\rm in}\,\,\, B_{2R} \\
 u_\ep =u  &
{\rm on}\,\,\,
\partial B_{2R} \,.
\end{cases}
\end{equation}
Since $u \in C^{1,\alpha}(\overline {B_{2R}})$, by \cite[Theorem 1.7 and subsequent remarks]{Li}
\begin{equation}\label{loc16}
\|u_\ep\|_{C^{1,\beta}(\overline {B_{2R}})} \leq C
\end{equation}
for some constant independent of $\varepsilon$. Hence, in particular, 
\begin{equation}\label{loc17}
\|a(|\nabla u_\ep|)\nabla u_\ep\|_{L^1(B_{2R})}\leq C
\end{equation}
for some constant independent of $\varepsilon$. The functions $a_\ep$ satisfy the assumptions imposed on $a$ in Step 1. Thus, by inequality \eqref{loc15},
\begin{align}\label{loc18}
\|a_\ep(|\nabla u_\ep|)\nabla u_\ep\|_{W^{1,2}(B_R)} \leq C\big(\|f\|_{B_{2R}} + R^{-\frac n2}\|a_\ep(|\nabla u_\ep|)\nabla u_\ep\|_{L^1(B_{2R})}\big)\,,
\end{align}
where, owing to \eqref{indici}, the constant  $C=C(n, i_a, s_a)$, and, in particular, is indepedent of $\ep$. Inequalities \eqref{loc17} and \eqref{loc18} ensure that the sequence $\{a_\ep(|\nabla u_\ep|)\nabla u_\ep\}$ is bounded in $W^{1,2}(B_R)$, and hence there exists a function $U: B_R \to \mathbb R^n$, with $U\in W^{1,2}(B_R)$, and a sequence $\{\ep _k\}$ such that
\begin{equation}\label{loc19}
a_{\ep _k}(|\nabla u_{\ep _k}|) \nabla u_{\ep _k} \to U\quad \hbox{in $L^2(B_R)$} \quad \hbox{and} \quad a_{\ep _k}(|\nabla u_{\ep _k}|) \nabla u_{\ep _k} \rightharpoonup U\quad \hbox{in $W^{1,2}(B_R)$}.
\end{equation}
Moreover, by \eqref{loc16}, there exists a function $v \in C^1(\overline {B_{2R}})$ such that, up to  subsequences,  
\begin{equation}\label{loc21}
 u_{\ep _k} \to v \quad \hbox{and} \quad \nabla u_{\ep _k} \to \nabla v
\end{equation}
pointwise in $\overline {B_{2R}}$.  In particular, 
\begin{equation}\label{loc21bis}
v=u \quad \hbox{on $\partial {B_{2R}}$,}
\end{equation}
inasmuch as $u_{\ep _k} =u$ on $\partial {B_{2R}}$ for every $k \in \mathbb R$.
Thanks to \eqref{loc19} and \eqref{loc21},
\begin{equation}\label{loc22}
a(|\nabla v|) \nabla v = U \in W^{1,2}(B_R).
\end{equation}
The weak formulation of problem \eqref{eqeploc} amounts to
\begin{equation}\label{loc23}
\int _{B_{2R}} a_{\ep_k}(|\nabla u_{\ep_k}|) \nabla u_{\ep_k} \cdot \nabla \varphi \, dx = \int_{B_{2R}} f \, \varphi \, dx
\end{equation}
for every $\varphi \in C^\infty_0(B_{2R})$. By \eqref{loc19} and \eqref{loc22}, 
passing to the limit in \eqref{loc23} as $k \to \infty$ results in
%
\begin{equation}\label{loc24}
\int _{B_{2R}} a(|\nabla v|) \nabla v \cdot \nabla \varphi \, dx = \int_{B_{2R}} f \, \varphi \, dx\,.
\end{equation}
Thus $v$ is the weak solution to the problem 
\begin{equation}\label{loc25}
\begin{cases}
- {\rm div} (a(|\nabla v|)\nabla v ) = f(x)  & {\rm in}\,\,\, B_{2R} \\
 v =u  &
{\rm on}\,\,\,
\partial B_{2R} \,.
\end{cases}
\end{equation}
Since $u$ solves the same problem,  $u=v$ in $B_{2R}$. Moreover, equations \eqref{loc18}, \eqref{loc19} and   \eqref{loc22} entail that $a(|\nabla u|) \nabla u   \in W^{1,2}(B_R)$, and 
\begin{align}\label{loc26}
\|a(|\nabla u|)\nabla u\|_{W^{1,2}(B_R)} \leq C\big(\|f\|_{L^2(B_{2R})} + R^{-\frac n2}\|a(|\nabla u|)\nabla u\|_{L^1(B_{2R})}\big)\,.
\end{align}
\par\noindent
\emph{Step 3}. Let $a$ and $f$ be as in the statement,  let $u$ be a   generalized local solution  to equation \eqref{localeq}, and let $f_k$ and $u_k$ be as in the definition of this kind of solution given at the begining of the present section.
An application of Step 2 to $u_k$ tells us that 
$a(|\nabla u_k|) \nabla u_k   \in W^{1,2}(B_R)$, and 
\begin{align}\label{loc27}
\|a(|\nabla u_k|)\nabla u_k\|_{W^{1,2}(B_R)}&  \leq C\big(\|f_k\|_{L^2(B_{2R})} + R^{-\frac n2}\|a(|\nabla u_k|)\nabla u_k\|_{L^1(B_{2R})}\big)
\\ \nonumber & \leq  C\big(\|f_k\|_{L^2(B_{2R})} + R^{-\frac n2}\|a(|\nabla u|)\nabla u\|_{L^1(B_{2R})}\big)
\,,
\end{align}
where the constant $C$ is independent of $k$. Therefore, the sequence $\{a(|\nabla u_k|)\nabla u_k\}$ is bounded in $W^{1,2}(B_R)$, and hence 
there exists a function
$U: B_R \to \mathbb R^n$, with $U\in W^{1,2}(B_R)$, and a subsequence, still indexed by $k$,   such that
\begin{equation}\label{loc28}
a(|\nabla u_k|) \nabla u_k \to U\quad \hbox{in $L^2(B_R)$} \quad \hbox{and} \quad a_k(|\nabla u_k|) \nabla u_k \rightharpoonup U\quad \hbox{in $W^{1,2}(B_R)$}.
\end{equation}
By assumption \eqref{approxuloc}, $\nabla u_k \to \nabla u$ a.e. in $\Omega$. Hence, owing to \eqref{loc28}, 
\begin{equation}\label{loc30}
a(|\nabla u|) \nabla u = U \quad \hbox{in $B_R$,}
\end{equation}
and
\begin{equation}\label{loc31}
\liminf _{k \to \infty} \|a(|\nabla u_k|)\nabla u_k\|_{W^{1,2}(B_R)} \geq \|a(|\nabla u|)\nabla u\|_{W^{1,2}(B_R)}\,.
\end{equation}
Inequality \eqref{secondloc2} follows from \eqref{loc27} and \eqref{loc31}. \qed

\end{document}